# QUADRATIC GROWTH CONDITIONS FOR CONVEX MATRIX OPTIMIZATION PROBLEMS ASSOCIATED WITH SPECTRAL FUNCTIONS

YING CUI*, CHAO DING†, AND XINYUAN ZHAO‡

**Abstract.** In this paper, we provide two types of sufficient conditions for ensuring the quadratic growth conditions of a class of constrained convex symmetric and non-symmetric matrix optimization problems regularized by nonsmooth spectral functions. These sufficient conditions are derived via the study of the $\mathcal{C}^2$-cone reducibility of spectral functions and the metric subregularity of their subdifferentials, respectively. As an application, we demonstrate how quadratic growth conditions are used to guarantee the desirable fast convergence rates of the augmented Lagrangian methods (ALM) for solving convex matrix optimization problems. Numerical experiments on an easy-to-implement ALM applied to the fastest mixing Markov chain problem are also presented to illustrate the significance of the obtained results.

**1. Introduction.** The *quadratic growth condition* is an important concept in optimization. It is closely related to the *metric subregularity* and *calmness* of set-valued mappings (see Section 2 for definitions), and the existence of *error bounds*. From different perspectives, the study of the metric subregularity and the calmness of set-valued mappings plays a central role in variational analysis, such as nonsmooth calculus and perturbation analysis of variational problems. We refer the reader to the monograph by Dontchev and Rockafellar [17] for a comprehensive study on both theory and applications of related subjects. See also [31, 18, 24, 37, 25, 26, 43] and references therein for recent advances.

Instead of considering general set-valued mappings, in this paper, we mainly focus on the solution mappings of convex matrix optimization problems. It is known from [1] that for convex problems, the metric subregularity of the subdifferentials of the essential objective functions (or the calmness of the solution mappings) can be equivalently characterized by the corresponding quadratic growth conditions. This connection motivates us to study the sufficient conditions for ensuring the latter properties. Beyond their own interests in second order variational analysis, the quadratic growth conditions can be employed for deriving the convergence rates of various first order and second order algorithms, including the proximal gradient methods [39, 58], the proximal point algorithms [50, 40, 36, 37] and the generalized Newton-type methods [22, 19, 42].

The convex matrix optimization problems concerned in our paper take the fol-

*Department of Mathematics, National University of Singapore, 10 Lower Kent Ridge Road, Singapore (matcuiy@nus.edu.sg).

†Institute of Applied Mathematics, Academy of Mathematics and Systems Science, Chinese Academy of Sciences, Beijing, P.R. China (dingchao@amss.ac.cn). The research of this author was supported by the National Natural Science Foundation of China under projects No. 11671387 and No. 11531014.

‡Beijing Institute for Scientific and Engineering Computing, Beijing University of Technology, Beijing, P.R. China. (xyzhao@bjut.edu.cn). The research of this author was supported by the Youth Program of National Natural Science Foundation of China under project No. 11101016.

lowing form:

$$\text{(1)} \quad \min_{X \in \mathbb{X}} \quad \Phi(X) := h(\mathcal{F}X) + \langle C, X \rangle + \theta(X)$$
$$\text{s.t.} \quad \mathcal{A}X \in b + \mathcal{Q},$$

where $\mathbb{X}$ is either the rectangular matrix space $\mathbb{R}^{n \times m}$ ($n \leqslant m$) or the symmetric matrix space $\mathbb{S}^n$, $C \in \mathbb{X}$ and $b \in \mathbb{R}^e$ are given data, $\mathcal{F} : \mathbb{X} \to \mathbb{R}^d$ and $\mathcal{A} : \mathbb{X} \to \mathbb{R}^e$ are two linear operators, $h : \mathbb{R}^d \to (-\infty, +\infty]$ is convex and continuously differentiable on $\text{dom } h$, which is assumed to be a nonempty open convex set, $\theta : \mathbb{X} \to (-\infty, +\infty]$ is a proper closed convex function and $\mathcal{Q} \subseteq \mathbb{R}^e$ is a given convex polyhedral cone. The dual of problem (1), in its equivalent minimization form, is given by

$$\text{(2)} \quad \min_{(y,w,S) \in \mathbb{R}^e \times \mathbb{R}^d \times \mathbb{X}} \quad \Psi(y, w, S) := \delta_{\mathcal{Q}^*}(y) - \langle b, y \rangle + h^*(-w) + \theta^*(-S)$$
$$\text{s.t.} \quad \mathcal{A}^* y + \mathcal{F}^* w + S = C,$$

where $\delta_{\mathcal{Q}^*}(\cdot)$ is the indicator function over the dual cone $\mathcal{Q}^*$, $\mathcal{A}^*$ and $\mathcal{F}^*$ are the adjoint of the linear operators $\mathcal{A}$ and $\mathcal{F}$, respectively, and $h^*$ and $\theta^*$ are the corresponding conjugate functions of $h$ and $\theta$. Problems of form (1) constitute a large class of convex matrix optimization problems with extensive applications in many fields, such as matrix completion, rank minimization, graph theory and machine learning. The (possibly nonsmooth) function $\theta$ in the objective can be used for imposing different properties to the decision variable $X$. Frequently used examples of $\theta$ include the indicator function over the positive semidefinite (PSD) cone in semidefinite programming [59], the nuclear norm function (i.e., the sum of all singular values of a matrix) in matrix completion problems [8, 9, 45], the spectral norm function (i.e., the largest singular value of a matrix) in matrix approximation problems [61, 28, 57], and the matrix Ky Fan $k$-norm ($1 \leqslant k \leqslant n$) function (i.e., the sum of $k$ largest singular values of a matrix) in fastest mixing Markov chain problems [7, 6].

The functions $\theta$ in the aforementioned examples belong to a special class of functions, called *spectral functions* [33, 35]. Specifically, they can be formulated either in the form of

$$\text{(3)} \quad \theta(X) = g(\sigma(X)), \quad X \in \mathbb{X} \ (\mathbb{R}^{n \times m} \text{ with } n \leqslant m \text{ or } \mathbb{S}^n)$$

with the function $g : \mathbb{R}^n \to (-\infty, +\infty]$ being proper closed convex and absolutely symmetric, or in the form of

$$\text{(4)} \quad \theta(X) = g(\lambda(X)), \quad X \in \mathbb{S}^n$$

with the function $g : \mathbb{R}^n \to (-\infty, +\infty]$ being proper closed convex and symmetric. Here $\sigma(X)$ denotes the vector of singular values for any given $X \in \mathbb{X}$ with the components $\sigma_1(X) \geqslant \sigma_2(X) \geqslant \ldots \geqslant \sigma_n(X) \geqslant 0$ being arranged in the non-increasing order, and $\lambda(X)$ denotes the vector of eigenvalues for any given $X \in \mathbb{S}^n$ with the components $\lambda_1(X) \geqslant \lambda_2(X) \geqslant \ldots \geqslant \lambda_n(X)$ also being arranged in the non-increasing order. In particular, for the indicator function over the PSD cone $\theta(\cdot) = \delta_{\mathbb{S}^n_+}(\cdot)$, the corresponding symmetric function $g$ is the indicator function over the positive orthant $\mathbb{R}^n_+$, and for the matrix Ky Fan $k$-norm function $\theta(\cdot) = \|\cdot\|_{(k)}$, the corresponding absolutely symmetric function $g$ is the sum of $k$ largest absolute components of a given vector.

In this paper, we concentrate on the study of sufficient conditions for guaranteeing the quadratic growth conditions for problem (1) and problem (2) associated with spectral functions (3) or (4). The sufficient conditions obtained in this work are of two types. One is based on the "no-gap" second-order sufficient conditions [5] of problems (1) and (2), under the $\mathcal{C}^2$-cone reducibility assumption of spectral functions. The other is through the bounded linear regularity of a collection of sets, whose intersection exactly gives the optimal solutions, under the metric subregularity assumption of the subdifferential mappings of spectral functions. The latter type has been used in [64] for studying the error bound of unconstrained convex matrix optimization problems, and in [12] for discussing the metric subregularity of the set-valued mappings arising from constrained semidefinite programming. Moreover, we show that the $\mathcal{C}^2$-cone reducibility of spectral functions, as well as the metric subregularity of their subdifferential mappings, can be justified by the corresponding properties of underlying (absolutely) symmetric functions. These results are also of their own interests as they provide tools for verifying these two fundamental variational properties of the non-polyhedral spectral function $\theta$ through the possibly polyhedral (absolutely) symmetric function $g$. In particular, for all the aforementioned examples, these two properties with respect to $g$ hold automatically due to their piecewise linear structures.

To illustrate the usefulness of our derived results, we investigate the fast convergence rates of the augmented Lagrangian method (ALM) for solving problem (2) under the quadratic growth conditions. Our motivation of this part stemmed from the highly promising numerical results of the ALM incorporated with the semismooth Newton-CG algorithm for solving large scale convex matrix problems [63, 32, 10, 62, 38]. We extend the results in the current literatures on the rates of the ALM for solving convex optimization problems [49, 40] and show that the (super)linear convergence rates of the ALM may still be valid even if problem (1) admits multiple solutions.

The remaining parts of this paper are organized as follows. In the next section, we introduce some preliminary knowledge from variational analysis in formulations and proofs of the main results. In Section 3, we establish the quadratic growth conditions of problem (1) (or problem (2)) under the assumptions that either $g$ (or $g^*$) is $\mathcal{C}^2$-cone reducible or $\partial g$ (or $\partial g^*$) is metrically subregular. Section 4 is devoted to an application of the quadratic growth conditions for the convex matrix optimization problems, that is, we establish the asymptotic (super)linear convergence rates of the ALM under the quadratic growth conditions. In Section 5, we conduct numerical experiments on solving fastest mixing Markov chain problems to demonstrate the derived fast rates of the ALM. We conclude our paper and make some comments in the final section.

**2. Notation and preliminary.** Let $\mathbb{U}$ and $\mathbb{V}$ be two finite dimensional real Euclidean spaces. For any $u \in \mathbb{U}$ and $\rho > 0$, we define the ball $\mathbb{B}_\rho(u) := \{v \in \mathbb{U} \mid \|v - u\| \leq \rho\}$. Let $D \subseteq \mathbb{U}$ be a set. For any $u \in D$, the tangent cone of $D$ at $u$ is defined by $\mathcal{T}_D(u) := \{d \in \mathbb{U} \mid \exists u^k \to u \text{ with } u^k \in D \text{ and } t^k \downarrow 0 \text{ such that } (u^k - u)/t^k \to d\}$. We let $\delta_D(\cdot)$ to be the indicator function over $D$, i.e., $\delta_D(x) = 0$ if $x \in D$, and $\delta_D(x) = \infty$ if $x \notin D$. If $D \subseteq \mathbb{U}$ is a convex set, we use $\text{ri}(D)$ to denote its relative interior. For a given closed convex set $D \subseteq \mathbb{U}$ and $u \in \mathbb{U}$, define $\Pi_D(u) := \arg\min\{\|d - u\| \mid d \in D\}$ and $\text{dist}(u, D) := \min_{d \in D} \|d - u\|$. Let $\alpha \subseteq \{1, ..., n\}$ and $\beta \subseteq \{1, ..., m\}$ be two index sets. For any $Z \in \mathbb{R}^{n \times m}$, we write $Z_i$ to be the $i$-th column of $Z$ and $Z_{\alpha\beta}$ to be the $|\alpha| \times |\beta|$ sub-matrix of $Z$ obtained by removing all the rows of $Z$ not in $\alpha$ and all the columns of $Z$ not in $\beta$. For any $z \in \mathbb{R}^n$, we denote $\text{Diag}(z)$ as the $n \times n$ diagonal matrix whose $i$-th diagonal entry is $z_i$ for $i = 1, 2, \ldots, n$. For any $Z \in \mathbb{R}^{n \times n}$, we denote $\text{diag}(Z)$ as the column vector consisting of all the diagonal entries of $Z$ being arranged from the first to the last. For a given proper closed convex function $p : \mathbb{U} \to (-\infty, +\infty]$, we use $\text{dom}\, p$ to denote its effective domain, $\text{epi}\, p$ to denote its epigraph, $p^*$ to denote its conjugate and $\partial p$ to denote its subdifferential, as in standard convex analysis [48]. We also use $\text{Prox}_p$ to denote the proximal mapping of $p$.

For a given positive integer $n$, let $\mathbb{O}^n$ be the set of all $n \times n$ orthogonal matrices. For any $X \in \mathbb{X}$, let $\mathbb{O}_{\mathbb{X}}(X)$ be the set of paired orthogonal matrices satisfying the singular value decomposition, i.e.,

$$\mathbb{O}_{\mathbb{X}}(X) = \begin{cases} \{(U,V) \in \mathbb{O}^n \times \mathbb{O}^m \mid X = U[\text{Diag}(\sigma(X))\ 0]V^T\} & \text{if } \mathbb{X} = \mathbb{R}^{n \times m}, \\ \{(U,V) \in \mathbb{O}^n \times \mathbb{O}^n \mid X = U\text{Diag}(\sigma(X))V^T\} & \text{if } \mathbb{X} = \mathbb{S}^n. \end{cases}$$

To distinguish with the above sets of singular vectors, we use $\mathcal{O}_{\mathbb{S}^n}(X)$ to denote the set of orthogonal matrices satisfying the eigenvalue decomposition of a given matrix $X \in \mathbb{S}^n$, i.e.,

$$\mathcal{O}_{\mathbb{S}^n}(X) = \{P \in \mathbb{O}^n \mid X = P\text{Diag}(\lambda(X))P^T\}.$$

The following lemma regarding the upper Lipschitz continuity of singular vectors and eigenvectors is given in [16, Proposition 7] and [55, Lemma 4.3], respectively.

LEMMA 1. *(i) Let $\overline{X} \in \mathbb{X}$. Then there exist constants $\varepsilon > 0$ and $\kappa > 0$ such that for any $X \in \mathbb{B}_\varepsilon(\overline{X})$ and any $(U,V) \in \mathbb{O}_{\mathbb{X}}(X)$, there exists $(\overline{U},\overline{V}) \in \mathbb{O}_{\mathbb{X}}(\overline{X})$ such that*

$$\|(U,V) - (\overline{U},\overline{V})\| \leq \kappa \|X - \overline{X}\|.$$

*(ii) Let $\overline{X} \in \mathbb{S}^n$. Then there exist constants $\varepsilon > 0$ and $\kappa > 0$ such that for any $X \in \mathbb{B}_\varepsilon(\overline{X})$ and any $P \in \mathcal{O}_{\mathbb{S}^n}(X)$, there exists $\overline{P} \in \mathcal{O}_{\mathbb{S}^n}(\overline{X})$ such that*

$$\|P - \overline{P}\| \leq \kappa \|X - \overline{X}\|.$$

A function $g : \mathbb{R}^n \to (-\infty, +\infty]$ is said to be *symmetric* if $g(x) = g(Qx)$ for any $x \in \mathbb{R}^n$ and any permutation matrix $Q \in \mathbb{R}^{n \times n}$, and is said to be *absolutely symmetric* if $g(x) = g(Qx)$ for any $x \in \mathbb{R}^n$ and any signed permutation matrix $Q \in \mathbb{R}^{n \times n}$ (i.e., an $n \times n$ matrix each of whose rows and columns has one nonzero entry which is $\pm 1$). The following two propositions, which are taken from [33, 34, 35], provide the formulas for the conjugate and the subdifferentials of spectral functions. In particular, part (ii) in Proposition 2 generalizes the characterization of the subdifferentials for the orthogonally invariant norms in [60].

PROPOSITION 2. *Let $g : \mathbb{R}^n \to (-\infty, +\infty]$ be a proper closed convex and absolutely symmetric function.*
*(i) The conjugate function $g^*$ is absolutely symmetric and $(g \circ \sigma)^* = g^* \circ \sigma$.*
*(ii) Let $X \in \mathbb{X}$ have the singular value $\sigma(X)$ in $\text{dom}\, g$. Let $W \in \mathbb{X}$. Then $W \in \partial(g \circ \sigma)(X)$ if and only if $\sigma(W) \in \partial g(\sigma(X))$ and there exists $(U,V) \in \mathbb{O}_{\mathbb{X}}(X) \cap \mathbb{O}_{\mathbb{X}}(W)$. In fact, $\partial(g \circ \sigma)(X) = \{U[\text{Diag}(\mu)\ 0]V^T \mid \mu \in \partial g(\sigma(X)), (U,V) \in \mathbb{O}_{\mathbb{X}}(X)\}$.*

PROPOSITION 3. *Let $g : \mathbb{R}^n \to (-\infty, +\infty]$ be a proper closed convex and symmetric function.*
*(i) The conjugate function $g^*$ is symmetric and $(g \circ \lambda)^* = g^* \circ \lambda$.*
*(ii) Let $X \in \mathbb{S}^n$ have the eigenvalue $\lambda(X)$ in $\text{dom}\, g$. Let $W \in \mathbb{S}^n$. Then $W \in \partial(g \circ \lambda)(X)$ if and only if $\lambda(W) \in \partial g(\lambda(X))$ and there exists $P \in \mathcal{O}_{\mathbb{S}^n}(X) \cap \mathcal{O}_{\mathbb{S}^n}(W)$. In fact, $\partial(g \circ \lambda)(X) = \{P\text{Diag}(\mu)P^T \mid \mu \in \partial g(\lambda(X)), P \in \mathcal{O}_{\mathbb{S}^n}(X)\}$.*

Let $G : \mathbb{U} \rightrightarrows \mathbb{V}$ be a set-valued mapping. The graph of $G$ is defined as $\text{gph}\, G := \{(u,v) \in \mathbb{U} \times \mathbb{V} \mid v \in G(u)\}$ and the inverse mapping of $G$ is defined as $G^{-1}(v) = \{u \in \mathbb{U} \mid v \in G(u)\}$ for any $v \in \mathbb{V}$. The following definition of metric subregularity is taken from [17, Section 3.8(3H)].

DEFINITION 4. *Let $\mathbb{U}$ and $\mathbb{V}$ be two finite dimensional Euclidean spaces. A set-valued mapping $G : \mathbb{U} \rightrightarrows \mathbb{V}$ is called metrically subregular at $\bar{u}$ for $\bar{v}$ (with modulus $\kappa$) if $(\bar{u}, \bar{v}) \in \mathrm{gph}\, G$ and there exist constants $\delta > 0$, $\varepsilon > 0$ and $\kappa > 0$ such that*

$$\mathrm{dist}(u, G^{-1}(\bar{v})) \leqslant \kappa\, \mathrm{dist}(\bar{v}, G(u) \cap \mathbb{B}_\delta(\bar{v})) \quad \forall\, u \in \mathbb{B}_\varepsilon(\bar{u}).$$

It is known that (see, e.g., [17, Theorem 3H.3]) for any set-valued mapping $G : \mathbb{U} \rightrightarrows \mathbb{V}$ and any $(\bar{u}, \bar{v}) \in \mathrm{gph}\, G$, the mapping $G$ is metrically subregular at $\bar{u}$ for $\bar{v}$ if and only if $G^{-1}$ is calm at $\bar{v}$ for $\bar{u}$, i.e., there exist constants $\delta' > 0$, $\varepsilon' > 0$ and $\kappa' > 0$ such that

$$G^{-1}(v) \cap \mathbb{B}_{\delta'}(\bar{u}) \subseteq G^{-1}(\bar{v}) + \kappa'\|v - \bar{v}\|\mathbb{B}_\mathbb{U} \quad \forall\, v \in \mathbb{B}_{\varepsilon'}(\bar{v}),$$

where $\mathbb{B}_\mathbb{U}$ is the unit ball in $\mathbb{U}$.

The equivalence between the quadratic growth condition of a proper closed convex function and the metric subregularity of its subdifferential has first been proved in [1, Theorem 3.3] in Hilbert spaces and then been extended in [2, Theorem 2.1] to Banach spaces. We will restrict our attention to this equivalence in Euclidean spaces.

PROPOSITION 5. *Suppose that $p : \mathbb{U} \to (-\infty, +\infty]$ is a proper closed convex function. Let $(\bar{x}, \bar{v}) \in \mathbb{U} \times \mathbb{U}$ satisfy $(\bar{x}, \bar{v}) \in \mathrm{gph}\, \partial p$. Then $\partial p$ is metrically subregular at $\bar{x}$ for $\bar{v}$ if and only if there exist constants $\kappa > 0$ and $\delta > 0$ such that*

$$(5) \qquad p(x) \geqslant p(\bar{x}) + \langle \bar{v}, x - \bar{x}\rangle + \kappa\, \mathrm{dist}^2(x, (\partial p)^{-1}(\bar{v})) \quad \forall\, x \in \mathbb{B}_\delta(\bar{x}).$$

*Specifically, if (5) holds with constant $\kappa$, then $\partial p$ is metrically subregular at $\bar{x}$ for $\bar{v}$ with modulus $1/\kappa$; conversely, if $\partial p$ is metrically subregular at $\bar{x}$ for $\bar{v}$ with modulus $\kappa'$, then (5) holds for all $\kappa \in (0, 1/(4\kappa'))$.*

The following $\mathcal{C}^2$-cone reducible property is adopted from [5, Definition 3.135] and is needed in our subsequent discussions.

DEFINITION 6. *Let $\mathcal{Q} \subseteq \mathbb{U}$ be a pointed convex closed cone (a cone is said to be pointed if $z \in \mathcal{Q}$ and $-z \in \mathcal{Q}$ implies that $z = 0$). The closed convex set $\mathcal{K} \subseteq \mathbb{V}$ is said to be $\mathcal{C}^2$-cone reducible at $\overline{X} \in \mathcal{K}$ to the cone $\mathcal{Q}$, if there exist an open neighborhood $\mathcal{W} \subseteq \mathbb{V}$ of $\overline{X}$ and a twice continuously differentiable mapping $\Xi : \mathcal{W} \to \mathbb{U}$ such that: (i) $\Xi(\overline{X}) = 0 \in \mathbb{U}$; (ii) the derivative mapping $\Xi'(\overline{X}) : \mathbb{V} \to \mathbb{U}$ is onto; (iii) $\mathcal{K} \cap \mathcal{W} = \{X \in \mathcal{W} \mid \Xi(X) \in \mathcal{Q}\}$. We say that $\mathcal{K}$ is $\mathcal{C}^2$-cone reducible if $\mathcal{K}$ is $\mathcal{C}^2$-cone reducible at every $\overline{X} \in \mathcal{K}$.*

Based on the above definition, we say a proper closed convex function $p : \mathbb{U} \to (-\infty, \infty]$ is $\mathcal{C}^2$-cone reducible at $u \in \mathrm{dom}\, p$ if $\mathrm{epi}\, p$ is $\mathcal{C}^2$-cone reducible at $(u, p(u))$. Moreover, $p$ is said to be $\mathcal{C}^2$-cone reducible if it is $\mathcal{C}^2$-cone reducible at every $u \in \mathrm{dom}\, p$.

We recall the following definition of bounded linear regularity (see, e.g., [3]).

DEFINITION 7. *Let $D_1, D_2, \ldots, D_s \subseteq \mathbb{U}$ be closed convex sets for some positive integer $s$. Suppose that $D := D_1 \cap D_2 \cap \ldots \cap D_s$ is non-empty. The collection $\{D_1, D_2, \ldots, D_s\}$ is said to be boundedly linearly regular if for every bounded set $\mathcal{B} \subseteq \mathbb{U}$, there exists a constant $\kappa > 0$ such that*

$$\mathrm{dist}(x, D) \leqslant \kappa \max\{\mathrm{dist}(x, D_1), \mathrm{dist}(x, D_2), \ldots, \mathrm{dist}(x, D_s)\} \quad \forall\, x \in \mathcal{B}.$$

The property of bounded linear regularity can be implied by the standard constraint qualification, as stated in the following proposition [4, Corollary 3].

PROPOSITION 8. *Let $D_1, D_2, \ldots, D_s \subseteq \mathbb{U}$ be closed convex sets for some positive integer $s$. Suppose that $D_1, D_2, \ldots, D_{s_0}$ are polyhedral sets for some $s_0 \in \{0, 1, \ldots, s\}$. Then a sufficient condition for $\{D_1, D_2, \ldots, D_s\}$ to be boundedly linearly regular is*

$$\bigcap_{i=1,2,\ldots,s_0} D_i \;\cap\; \bigcap_{i=s_0+1,\ldots,s} \mathrm{ri}(D_i) \neq \emptyset.$$

Finally, for any function $p : \mathbb{U} \to (-\infty, +\infty]$, we use $p_{-}^{\downarrow}(u; d)$ and $p_{+}^{\downarrow}(u; d)$ to denote the lower and upper directional epiderivatives of $p$ at $u \in \mathbb{U}$ along the direction $d \in \mathbb{U}$, i.e.,

$$p_{-}^{\downarrow}(u; d) = \liminf_{\substack{t \downarrow 0 \\ d' \to d}} \frac{p(u + td') - p(u)}{t}$$

and

$$p_{+}^{\downarrow}(u; d) = \sup_{\{t_n\} \in \Delta} \left\{ \liminf_{\substack{t_n \downarrow 0 \\ d' \to d}} \frac{p(u + t_n d') - p(u)}{t_n} \right\},$$

where $\Delta$ denotes the set of all positive real sequences $\{\rho_n\}$ converging to 0 [5, Section 2.2.3]. If $p$ is a closed convex function, then we know from [5, Proposition 2.58] that $p_{-}^{\downarrow}(u; \cdot) = p_{+}^{\downarrow}(u; \cdot)$ for any $u \in \text{dom}\, p$, i.e., $p$ is directionally epidifferentiable with the directional epiderivative $p^{\downarrow}(u; \cdot)$. Furthermore, if $p^{\downarrow}(u; d)$ is finite for $u \in \text{dom}\, p$ and $d \in \mathbb{U}$, we define the lower second order directional epiderivative of $p$ by for any $w \in \mathbb{U}$,

$$p_{-}^{\downarrow\downarrow}(u; d, w) := \liminf_{\substack{t \downarrow 0 \\ w' \to w}} \frac{p(u + td + \frac{1}{2}t^2 w') - p(u) - t p^{\downarrow}(u; d)}{\frac{1}{2}t^2}.$$

**3. The quadratic growth conditions for convex matrix optimization problems.** For notational simplicity, we write $\mathbb{Z} := \mathbb{R}^e \times \mathbb{R}^d \times \mathbb{X}$ and for any $y \in \mathbb{R}^e$, $w \in \mathbb{R}^d$ and $S \in \mathbb{X}$, we write $Z = (y, w, S)$. The Lagrangian function $l$ associated with problem (2) takes the form of

$$l(Z, X) := \Psi(Z) + \langle X, \mathcal{A}^* y + \mathcal{F}^* w + S - C \rangle, \quad (Z, X) \in \mathbb{Z} \times \mathbb{X}.$$

Let $\phi : \mathbb{X} \to (-\infty, \infty]$ be the essential objective function of problem (1), i.e.,

(6) $\quad \phi(X) := -\inf_{Z \in \mathbb{Z}} l(Z, X) = \begin{cases} \Phi(X) & \text{if } \mathcal{A} X \in b + \mathcal{Q}, \\ +\infty & \text{otherwise,} \end{cases} \quad \forall X \in \mathbb{X}.$

Denote the set-valued mapping $\mathcal{T}_\phi : \mathbb{X} \rightrightarrows \mathbb{X}$ associated with $\phi$ as

(7) $\quad \mathcal{T}_\phi(X) := \partial \phi(X), \quad X \in \mathbb{X}.$

Let $\Omega_P \subseteq \mathbb{X}$ and $\Omega_D \subseteq \mathbb{Z}$ be the optimal solution sets of the primal problem (1) and the dual problem (2), respectively, both being assumed to be nonempty. Let $\mathcal{M}_D(\overline{Z}) \subseteq \mathbb{X}$ be the set of Lagrangian multipliers associated with $\overline{Z} \in \Omega_D$ for problem (2), i.e., $\overline{X} \in \mathcal{M}_D(\overline{Z})$ if and only if $(\overline{X}, \overline{Z}) = (\overline{X}, \bar{y}, \bar{w}, \overline{S})$ solves the following KKT system:

(8) $\quad \begin{cases} 0 \in \mathcal{A} X - b + \mathcal{N}_{\mathcal{Q}^*}(y), \ 0 \in \mathcal{F} X - \partial h^*(-w), \ 0 \in X - \partial \theta^*(-S), \\ 0 = C - (\mathcal{A}^* y + \mathcal{F}^* w + S), \end{cases}$

for $(Z, X) \in \mathbb{Z} \times \mathbb{X}$, where $\mathcal{N}_{\mathcal{Q}^*}(y)$ denotes the normal cone of $\mathcal{Q}^*$ at $y \in \mathcal{Q}^*$ as in standard convex analysis [48]. It can be easily checked that if $(\overline{X}, \overline{Z}) \in \mathbb{X} \times \mathbb{Z}$ is a solution to the above KKT system, then $(\overline{X}, \bar{y})$ solves the following inclusion problem:

(9) $\quad \begin{cases} 0 \in C - \mathcal{A}^* y + \mathcal{F}^* \nabla h(\mathcal{F} X) + \partial \theta(X), \\ 0 \in \mathcal{A} X - b + \mathcal{N}_{\mathcal{Q}^*}(y), \end{cases} \quad (X, y) \in \mathbb{X} \times \mathbb{R}^e.$

We write $\mathcal{M}_P(\overline{X}) \subseteq \mathbb{R}^e$ as the set of Lagrangian multipliers $\bar{y}$ associated with $\overline{X} \in \Omega_P$, i.e., $\bar{y} \in \mathcal{M}_P(\overline{X})$ if and only if $(\overline{X}, \bar{y})$ satisfies (9).

Let $F_P$ and $F_D$ be the sets of all feasible points for problem (1) and problem (2), respectively, i.e.,

$$F_P = \{X \in \mathbb{X} \mid \mathcal{A}X \in b + \mathcal{Q}\}, \quad F_D = \{Z \in \mathbb{Z} \mid \mathcal{A}^*y + \mathcal{F}^*w + S = C, \; y \in \mathcal{Q}^*\}.$$

The quadratic growth condition is said to be satisfied at an optimal solution $\overline{X} \in \Omega_P$ for problem (1) if there exist constants $\delta_p > 0$ and $\kappa_p > 0$ such that

(10) $$\Phi(X) \geq \Phi(\overline{X}) + \kappa_p \operatorname{dist}^2(X, \Omega_P) \quad \forall\, X \in F_P \cap \mathbb{B}_{\delta_p}(\overline{X}).$$

Similarly, the quadratic growth condition is said to be satisfied at $\overline{Z} \in \Omega_D$ for problem (2) if there exist constants $\delta_d > 0$ and $\kappa_d > 0$ such that

(11) $$\Psi(Z) \geq \Psi(\overline{Z}) + \kappa_d \operatorname{dist}^2(Z, \Omega_D) \quad \forall\, Z \in F_D \cap \mathbb{B}_{\delta_d}(\overline{Z}).$$

The constants $\kappa_p$ in (10) and $\kappa_d$ in (11) are called the quadratic growth modulus with respect to the problem (1) and (2), respectively.

The lemma below, which is a direct application of Proposition 5 to the essential objective function $\phi$ in (6), reveals the relationship between the metric subregularity of $\mathcal{T}_\phi$ in (7) at an optimal solution for the origin and the quadratic growth condition for problem (1). An analogous result regarding the dual problem (2) can be established in the same fashion.

LEMMA 9. *The following two properties are equivalent to each other:*
*(a) The quadratic growth condition (10) holds at $\overline{X} \in \Omega_P$ for problem (1).*
*(b) The set-valued mapping $\mathcal{T}_\phi$ is metrically subregular at $\overline{X} \in \Omega_P$ for the origin. Specifically, if (10) holds with quadratic growth modulus $\kappa_p$, then $\mathcal{T}_\phi$ is metrically subregular at $\overline{X}$ for the origin with modulus $1/\kappa_p$; conversely, if $\mathcal{T}_\phi$ is metrically subregular at $\overline{X}$ for the origin with modulus $\kappa'_p$, then (10) holds for any $\kappa_p \in (0, 1/(4\kappa'_p))$.*

In the following two subsections, we shall study two types of sufficient conditions for ensuring the primal quadratic growth condition (10) and the dual quadratic growth condition (11), respectively. One is under the $\mathcal{C}^2$-cone reducibility assumption of $g$ ($g^*$), the other is under the metric subregularity assumption of $\partial g$ ($\partial g^*$).

**3.1. Sufficient conditions under the $\mathcal{C}^2$-cone reducibility.** When the spectral function $\theta$ in problem (1) is $\mathcal{C}^2$-cone reducible, one can derive the conditions for guaranteeing the quadratic growth condition (10) by employing the second order sufficient conditions of convex composite optimizations (see [5, Section 3.4.1] for details). In general, it may be difficult to check the $\mathcal{C}^2$-cone reducibility of matrix functions. Our first result is to show that the $\mathcal{C}^2$-cone reducibility of the spectral function $\theta$ can be verified via its underlying (absolutely) symmetric function $g$.

PROPOSITION 10. *Let $g : \mathbb{R}^n \to (-\infty, +\infty]$ be a proper closed convex and absolutely symmetric function. Let $\theta : \mathbb{X} \to (-\infty, +\infty]$ be the spectral function associated with $g$ as in (3). For any $\overline{X} \in \operatorname{dom}\theta$, if the function $g$ is $\mathcal{C}^2$-cone reducible at $\sigma(\overline{X})$, then the function $\theta$ is $\mathcal{C}^2$-cone reducible at $\overline{X}$.*

*Proof.* We first introduce some notations. Denote the nonzero distinct singular values of $\overline{X}$ by $\bar{\nu}_1 > \ldots > \bar{\nu}_r > 0$, where $r$ is a positive integer. Define the index sets $a_l := \{i \mid \sigma_i(\overline{X}) = \bar{\nu}_l, \; 1 \leq i \leq n\}$ for $l = 1, \ldots, r$, $a_{r+1} = b := \{i \mid \sigma_i(\overline{X}) = 0, \; 1 \leq i \leq n\}$ and $c := \{n+1, \ldots, m\}$. We write $\Sigma(X) = \operatorname{Diag}(\sigma(X))$ for any $X \in \mathbb{X}$ and $\overline{\Sigma} = \operatorname{Diag}(\sigma(\overline{X}))$.

Let $l \in \{1, \ldots, r+1\}$ be arbitrary but fixed. Since the singular value function $\sigma(\cdot)$ is globally Lipschitz continuous, one can obtain from [16, Proposition 5] that there exists an open neighborhood $\mathcal{N}$ of $\overline{X}$ such that the following functions $\mathcal{U}_l : \mathcal{N} \to \mathbb{R}^{n \times n}$ and $\mathcal{V}_l : \mathcal{N} \to \mathbb{R}^{m \times m}$ are well-defined:

$$(12) \quad \mathcal{U}_l(X) := \sum_{i \in a_l} U_i(X) U_i(X)^T, \quad \mathcal{V}_l(X) = \sum_{i \in a_l} V_i(X) V_i(X)^T, \quad X \in \mathcal{N},$$

where $(U(X), V(X)) \in \mathbb{O}_{\mathbb{X}}(X)$. That is, for any $X \in \mathcal{N}$, the function values $\mathcal{U}_l(X)$ and $\mathcal{V}_l(X)$ are independent of the choices of the orthogonal pairs $(U(X), V(X)) \in \mathbb{O}_{\mathbb{X}}(X)$. By the relationship between the singular value decomposition of any $X \in \mathbb{X}$ and the eigenvalue decomposition of its extended symmetric counterpart $\begin{bmatrix} 0 & X \\ X^T & 0 \end{bmatrix}$, one can derive from [16, Proposition 8] that the mappings $\mathcal{U}_l(\cdot)$ and $\mathcal{V}_l(\cdot)$ are twice continuously differentiable on $\mathcal{N}$.

To proceed, let us first consider the case that $\overline{X} = \begin{bmatrix} \Sigma & 0 \end{bmatrix}$. For any $X \in \mathcal{N}$, let $\mathcal{L}_l(X)$ and $\mathcal{R}_l(X)$ be the left and right singular vector spaces corresponding to the single values $\{\sigma_i(X) : i \in a_l\}$. Obviously the spaces spanned by the columns of $\mathcal{U}_l(X)$ and $\mathcal{V}_l(X)$ coincide with $\mathcal{L}_l(X)$ and $\mathcal{R}_l(X)$, respectively. Now we show that if $X$ is sufficiently close to $\overline{X}$, the columns of $\mathcal{U}_l(X)$ are the bases of $\mathcal{L}_l(X)$. In fact, for any $X \in \mathcal{N}$, the $i$-th column of $\mathcal{U}_l(X)$ is given by

$$(\mathcal{U}_l(X))_i = \sum_{j \in a_l} U_{ij}(X) \begin{bmatrix} U_{1j}(X) \\ \vdots \\ U_{nj}(X) \end{bmatrix}, \quad i \in a_l.$$

Let $\bar{q} = (\bar{q}_1, \ldots, \bar{q}_{|a_l|})^T \in \mathbb{R}^{|a_l|}$ be the vector such that $\sum_{i \in a_l} \bar{q}_i (\mathcal{U}_l(X))_i = 0$. Since the columns of $U(X)$ among the index set $a_l$ are linearly independent, we see that $\bar{q}$ is the solution of the linear system $U_{a_l a_l}(X) q = 0$. Recall that $\overline{X} = \begin{bmatrix} \Sigma & 0 \end{bmatrix}$. Then by [16, (31) in Proposition 7], for all $X$ sufficiently close to $\overline{X}$, there exists $Q_l \in \mathbb{O}^{|a_l|}$ such that $U_{a_l a_l}(X) = Q_l + O(\|X - \overline{X}\|)$. Thus, the matrix $U_{a_l a_l}(X)$ is invertible, which further implies that $\bar{q} = 0$ and the columns of $\mathcal{U}_l(X)$ are linearly independent. Consequently, by shrinking $\mathcal{N}$ if necessary, we have that for any $X \in \mathcal{N}$, the columns of $\mathcal{U}_l(X)$ are the bases of $\mathcal{L}_l(X)$. By using the similar arguments, one can also derive that the columns of $\mathcal{V}_l(X)$ are the bases of $\mathcal{R}_l(X)$. Applying the Gram-Schmidt orthonormalization procedure to the columns of $\mathcal{U}_l(X)$ and $\mathcal{V}_l(X)$ for $X \in \mathcal{N}$, we can obtain two matrices $M_{a_l}(X) \in \mathbb{R}^{n \times |a_l|}$ and $N_{a_l}(X) \in \mathbb{R}^{m \times |a_l|}$ such that the columns of $M_{a_l}(X)$ are the orthogonal bases of $\mathcal{L}_l(X)$ and the columns of $N_{a_l}(X)$ are the orthogonal bases of $\mathcal{R}_l(X)$. For $X$ sufficiently close to $\overline{X}$ and two positive integers $l, l' \in \{1, \ldots, r+1\}$ with $l \neq l'$, it holds that $\sigma_i(X) \neq \sigma_j(X)$ for $i \in a_l$ and $j \in a_{l'}$, implying that the matrices $M(x) = \begin{bmatrix} M_{a_1}(X) & \cdots & M_{a_r}(X) & M_{a_{r+1}}(X) \end{bmatrix} \in \mathbb{R}^{n \times n}$ and $N(x) = \begin{bmatrix} N_{a_1}(X) & \cdots & N_{a_r}(X) & N_{a_{r+1}}(X) \end{bmatrix} \in \mathbb{R}^{m \times m}$ are orthogonal and satisfy the following singular value decomposition

$$X = M(X) \begin{bmatrix} \Sigma(X) & 0 \end{bmatrix} N(X)^T.$$

Thus, we have

$$M_{a_l}(X)^T X N_{a_l}(X) = \begin{cases} (\Sigma(X))_{a_l a_l} & \text{for } l = 1, \ldots, r, \\ [\Sigma(X)_{bb} \quad 0] & \text{for } l = r+1. \end{cases}$$

Take $M_{a_l} : \mathcal{N} \to \mathbb{R}^{n \times |a_l|}$ and $N_{a_l} : \mathcal{N} \to \mathbb{R}^{n \times |a_l|}$ as two matrix mappings. Since $M_{a_l}$ and $N_{a_l}$ are twice continuously differentiable on $\mathcal{N}$, the mappings defined by $M_{a_l}(X)^T X N_{a_l}(X)$ for any $X \in \mathcal{N}$ are also twice continuously differentiable on $\mathcal{N}$. Moreover, one can derive from [16, (31) in Proposition 7] that for any $X \in \mathcal{N}$,

$$M_{a_l}(X) = \begin{bmatrix} O(\|X - \overline{X}\|) \\ I_{|a_l|} + O(\|X - \overline{X}\|^2) \\ O(\|X - \overline{X}\|) \end{bmatrix}, \quad N_{a_l}(X) = \begin{bmatrix} O(\|X - \overline{X}\|) \\ I_{|a_l|} + O(\|X - \overline{X}\|^2) \\ O(\|X - \overline{X}\|) \end{bmatrix}.$$

Denote $H := X - \overline{X}$. For any $X \in \mathcal{N}$, we deduce

$$(13) \quad M_{a_l}(X)^T X N_{a_l}(X) = \begin{cases} \overline{\Sigma}_{a_l a_l} + H_{a_l a_l} + O(\|H\|^2) & \text{for } l = 1, \ldots, r, \\ [\overline{\Sigma}_{bb} \; 0] + [H_{bb} \; H_{bc}] + O(\|H\|^2) & \text{for } l = r + 1. \end{cases}$$

For the general case, let $(\overline{U}, \overline{V}) \in \mathbb{O}_{\mathbb{X}}(\overline{X})$ be fixed. Then $\overline{U}^T X \overline{V} = [\overline{\Sigma} \; 0] + \overline{U}^T(X - \overline{X})\overline{V}$. Denote $\widetilde{H} = \overline{U}^T(X - \overline{X})\overline{V}$. Therefore, replacing $X$ by $\overline{U}^T X \overline{V}$ in the previous arguments, we know that there exists an open neighborhood $\mathcal{N}$ of $\overline{X}$ such that the mapping

$$\mathcal{D}_l(X) := \begin{cases} M_{a_l}(\overline{U}^T X \overline{V})^T \overline{U}^T X \overline{V} N_{a_l}(\overline{U}^T X \overline{V}) & \text{for } l = 1, \ldots, r, \\ M_{a_{r+1}}(\overline{U}^T X \overline{V})^T \overline{U}^T X \overline{V} N_{a_{r+1}}(\overline{U}^T X \overline{V}) & \text{for } l = r + 1, \end{cases} \quad X \in \mathcal{N},$$

$$= \begin{cases} (\Sigma(X))_{a_l a_l} & \text{for } l = 1, \ldots, r, \\ [\Sigma(X)_{bb} \; 0] & \text{for } l = r + 1, \end{cases}$$

is twice continuously differentiable on $\mathcal{N}$. In particular, we have

$$(14) \quad \mathcal{D}_l(\overline{X}) = \begin{cases} \overline{\Sigma}_{a_l a_l} & \text{for } l = 1, \ldots, r, \\ [\overline{\Sigma}_{bb} \; 0] = 0 & \text{for } l = r + 1. \end{cases}$$

Moreover, the equation (13) implies that for any $X \in \mathcal{N}$,

$$(15) \quad \mathcal{D}_l(X) - \mathcal{D}_l(\overline{X}) = \begin{cases} \widetilde{H}_{a_l a_l} + O(\|X - \overline{X}\|^2) & \text{for } l = 1, \ldots, r, \\ [\widetilde{H}_{bb} \; \widetilde{H}_{bc}] + O(\|X - \overline{X}\|^2) & \text{for } l = r + 1. \end{cases}$$

Let $\bar{t} = \theta(\overline{X})$ and $\mathcal{N}' \subseteq \mathbb{R}$ be any open neighborhood of $\bar{t}$ and $\mathcal{W} := \mathcal{N} \times \mathcal{N}'$. Define the mapping $\Xi : \mathcal{W} \to \mathbb{R}^n \times \mathbb{R}$ by

$$\Xi(X, t) := (\text{diag}(\mathcal{D}_1(X)), \ldots, \text{diag}(\mathcal{D}_r(X)), \text{diag}(\mathcal{D}_{r+1}(X)), t), \quad (X, t) \in \mathcal{N} \times \mathcal{N}'.$$

It then follows from (14) that $\text{epi}\,\theta \cap \mathcal{W} = \{(X, t) \in \mathcal{W} \mid \Xi(X, t) \in \text{epi}\,g\}$ and $\Xi(\overline{X}, \bar{t}) = (\sigma(\overline{X}), \bar{t})$. Moreover, (15) implies that for any $(H, \tau) \in \mathbb{X} \times \mathbb{R}$ sufficiently close to the origin,

$$\Xi(\overline{X} + H, \bar{t} + \tau) - \Xi(\overline{X}, \bar{t}) = (\text{diag}(\widetilde{H}_{a_1 a_1}), \ldots, \text{diag}(\widetilde{H}_{a_{r+1} a_{r+1}}), \tau) + O(\|(H, \tau)\|^2).$$

Thus, the derivative $\Xi'(\overline{X}, \bar{t}) : \mathbb{X} \times \mathbb{R} \to \mathbb{R}^n \times \mathbb{R}$ of $\Xi$ at $(\overline{X}, \bar{t})$ is given by

$$\Xi'(\overline{X}, \bar{t})(H, \tau) = (\text{diag}(\widetilde{H}_{a_1 a_1}), \ldots, \text{diag}(\widetilde{H}_{a_r a_r}), \text{diag}(\widetilde{H}_{bb}), \tau), \quad (H, \tau) \in \mathbb{X} \times \mathbb{R}.$$

Obviously $\Xi'(\overline{X}, \bar{t})$ is onto. Since $\text{epi}\,g$ is $\mathcal{C}^2$-cone reducible at $(\sigma(\overline{X}), \bar{t})$, we know from [52, Proposition 3.2] that $\text{epi}\,\theta$ is also $\mathcal{C}^2$-cone reducible at $(\overline{X}, \bar{t})$. The proof is thus completed. □

For spectral functions $\theta$ associated with symmetric functions $g$ in the form of (4), we have the following analogous conclusion. Its proof can be derived similarly to Proposition 10 by replacing all the singular value decompositions in the arguments with the eigenvalue decompositions. We omit the details here.

PROPOSITION 11. *Let $g : \mathbb{R}^n \to (-\infty, +\infty]$ be a proper closed convex and symmetric function. Let $\theta : \mathbb{S}^n \to (-\infty, +\infty]$ be the spectral function associated with $g$ as in (4). Then for any $\overline{X} \in \mathrm{dom}\, \theta$, if the function $g$ is $\mathcal{C}^2$-cone reducible at $\lambda(\overline{X})$, then the function $\theta$ is $\mathcal{C}^2$-cone reducible at $\overline{X}$.*

REMARK 1. *It is known that the convex polyhedral sets are $\mathcal{C}^2$-cone reducible [5, Example 3.139]. Hence, Propositions 10 and 11 imply that the class of $\mathcal{C}^2$-cone reducible spectral functions is rich. In particular, they cover the results in [5, Example 3.140] and [13, Proposition 4.3] about the $\mathcal{C}^2$-cone reducibility of the indicator function over the PSD cone and the matrix Ky Fan k-norm function, respectively.*

Let $\overline{X} \in \Omega_P$ with $\mathcal{M}_P(\overline{X}) \neq \varnothing$. The critical cone of problem (1) at $\overline{X}$ is given by
$$\mathcal{C}_P(\overline{X}) := \left\{ H \in \mathbb{X} \;\middle|\; \begin{array}{c} \mathcal{A}H \in \mathcal{T}_\mathcal{Q}(\mathcal{A}\overline{X} - b),\; H \in \mathcal{T}_{\mathrm{dom}\,\theta}(\overline{X}), \\ \langle \mathcal{F}^* \nabla h(\mathcal{F}\overline{X}) + C, H \rangle + \theta^\downarrow(\overline{X}; H) = 0 \end{array} \right\}.$$

Similarly, let $\overline{Z} \in \Omega_D$ with $\mathcal{M}_D(\overline{Z}) \neq \varnothing$. If $h^*$ is twice continuously differentiable, then the critical cone of problem (2) at $\overline{X}$ is given by
$$\mathcal{C}_D(\overline{Z}) := \left\{ \begin{array}{c} (H_y, H_w, H_S) \\ \in \mathbb{Z} \end{array} \;\middle|\; \begin{array}{c} \mathcal{A}^* H_y + \mathcal{F}^* H_w + H_S = 0, \\ H_y \in \mathcal{T}_{\mathcal{Q}^*}(\overline{y}),\; -H_S \in \mathcal{T}_{\mathrm{dom}\,\theta^*}(-\overline{S}) \\ \langle -b, H_y \rangle - \langle \nabla h^*(-\overline{w}), H_w \rangle + (\theta^*)^\downarrow(-\overline{S}; -H_S) = 0 \end{array} \right\}.$$

Now we are ready to present the main result of this subsection.

THEOREM 12. *Let $\theta$ be a spectral function in the form of (3) (or (4)). Let $\overline{X}$ be an optimal solution to problem (1) with $\mathcal{M}_P(\overline{X}) \neq \varnothing$. Assume that the following three conditions hold:*
*(a) the function $h$ is twice continuously differentiable on $\mathrm{dom}\, h$;*
*(b) the function $g$ is $\mathcal{C}^2$-cone reducible at $\sigma(\overline{X})$ (or $\lambda(\overline{X})$);*
*(c) the "no-gap" second order sufficient condition holds at $\overline{X}$ for problem (1), i.e., for any $H \in \mathcal{C}_P(\overline{X})\setminus\{0\}$,*
$$\sup_{\overline{y} \in \mathcal{M}_P(\overline{X})} \left\{ \langle \mathcal{F}H, \nabla^2 h(\mathcal{F}\overline{X}) \mathcal{F}H \rangle - \phi^*_{(\overline{X}, H)}(C - \mathcal{A}^* \overline{y} + \mathcal{F}^* \nabla h(\mathcal{F}\overline{X})) \right\} > 0,$$

*where $\phi^*_{(\overline{X}, H)}(\cdot)$ is the conjugate function of $\phi_{(\overline{X}, H)}(\cdot) := \theta^{\downarrow\downarrow}_-(\overline{X}; H, \cdot)$.*
*Then the quadratic growth condition (10) holds at the unique optimal solution $\overline{X}$ for problem (1).*

*Proof.* By applying Proposition 10 (Proposition 11) to condition (b), one can obtain the $\mathcal{C}^2$-cone reducibility of the spectral function $\theta$ at $\overline{X}$. Thus, we know from [5, Proposition 3.136] that $\mathrm{epi}\,\theta$ is second order regular at $(\overline{X}, \theta(\overline{X}))$ (see [5, Definition 3.85] for the definition of second order regular). Then the conclusion of Theorem 12 follows from [5, Theorem 3.109], directly. □

REMARK 2. *The explicit expressions of $\phi^*_{(\overline{X}, H)}(\cdot)$ in condition (c) for $\theta = \delta_{\mathbb{S}^n_+}(\cdot)$ can be found in [53] and for $\theta = \|\cdot\|_{(k)}$ can be found in [14]. Under the assumptions of Theorem 12, the point $\overline{X}$ is in fact the unique optimal solution of problem (1). Thus, the quadratic growth condition (10) at $\overline{X}$ can be equivalently stated as the existence of constants $\delta_p > p$ and $\kappa_p > 0$ such that*

(16) $$\Phi(X) \geqslant \Phi(\overline{X}) + \kappa_p \|X - \overline{X}\|^2 \quad \forall\, X \in F_P \cap \mathbb{B}_{\delta_p}(\overline{X}).$$

In the following, we also provide an analogous result to Theorem 12 regarding the dual quadratic growth condition (11). By taking into account Proposition 2(i) and Proposition 3(i), one can derive the proof in a similar manner.

THEOREM 13. *Let $\theta$ be a spectral function in the form of (3) (or (4)). Let $\overline{Z}$ be an optimal solution for problem (2) with $\mathcal{M}_D(\overline{Z}) \neq \emptyset$. Assume that the following three conditions hold:*
*(a) the function $h^*$ is twice continuously differentiable on dom $h^*$;*
*(b) the function $g^*$ is $\mathcal{C}^2$-cone reducible at $\sigma(-\overline{S})$ (or $\lambda(-\overline{S})$);*
*(c) the "no-gap" second order sufficient condition holds at $\overline{Z}$ for problem (2), i.e., for any $(H_y, H_w, H_S) \in \mathcal{C}_D(\overline{Z}) \backslash \{0\}$,*

$$(17) \qquad \sup_{\overline{X} \in \mathcal{M}_D(\overline{Z})} \left\{ \langle H_w, \nabla^2 h^*(-\bar{w}) H_w \rangle - \psi^*_{(\overline{S}, H_S)}(\overline{X}) \right\} > 0,$$

*where $\psi^*_{(\overline{S}, H_S)}(\cdot)$ is the conjugate function of $\psi_{(\overline{S}, H_S)}(\cdot) = (\theta^*)^{\downarrow\downarrow}_{-}(-\overline{S}; -H_S, \cdot)$.*
*Then the quadratic growth condition (11) holds at the unique optimal solution $\overline{Z}$ for problem (2).*

**3.2. Sufficient conditions under the metric subregularity of subdifferentials.** As discussed in Remark 2, the assumptions given in Theorem 12 (or Theorem 13) imply the uniqueness of the optimal solution for problem (1) (or problem (2)). In this subsection, we shall study the quadratic growth conditions from a different perspective. Our goal here is to provide sufficient conditions for this subject when the underlying problems admit multiple solutions. For this purpose, we first study the metric subregularity of the subdifferential of a spectral function $\theta$.

PROPOSITION 14. *Let $g : \mathbb{R}^n \to (-\infty, +\infty]$ be a proper closed convex and absolutely symmetric function. Let $\theta : \mathbb{X} \to (-\infty, +\infty]$ be the spectral function associated with $g$ as in (3). Consider any point $(\overline{X}, \overline{W}) \in \text{gph } \partial \theta$. If the subdifferential mapping $\partial g$ is metrically subregular at $\sigma(\overline{X})$ for $\sigma(\overline{W})$, then the subdifferential mapping $\partial \theta$ is metrically subregular at $\overline{X}$ for $\overline{W}$.*

*Proof.* Since $\partial g$ is assumed to be metrically subregular at $\sigma(\overline{X})$ for $\sigma(\overline{W})$, there exist constants $\delta > 0$, $\varepsilon > 0$ and $\kappa_0 > 0$ such that

$$(18) \quad \text{dist}(u, (\partial g)^{-1}(\sigma(\overline{W}))) \leqslant \kappa_0 \, \text{dist}(\sigma(\overline{W}), \partial g(u) \cap \mathbb{B}_\delta(\sigma(\overline{W}))) \quad \forall \, u \in \mathbb{B}_\varepsilon(\sigma(\overline{X})).$$

To establish the metric subregularity of $\partial \theta$ at $\overline{X}$ for $\overline{W}$, it suffices to show that there exists a constant $\kappa > 0$ such that

$$(19) \qquad \text{dist}(X, (\partial \theta)^{-1}(\overline{W})) \leqslant \kappa \, \text{dist}(\overline{W}, \partial \theta(X) \cap \mathbb{B}_\delta(\overline{W})) \quad \forall \, X \in \mathbb{B}_\varepsilon(\overline{X}).$$

Consider any $X \in \mathbb{B}_\varepsilon(\overline{X})$. If $\partial \theta(X) \cap \mathbb{B}_\delta(\overline{W}) = \emptyset$, the above inequality holds automatically. Thus, we only need to consider those $X$ such that $\partial \theta(X) \cap \mathbb{B}_\delta(\overline{W}) \neq \emptyset$. Let $W \in \partial \theta(X) \cap \mathbb{B}_\delta(\overline{W})$ be arbitrarily given. By Proposition 2 and the globally Lipschitz continuity of the singular value function $\sigma(\cdot)$, we see that

$$(20) \qquad \sigma(W) \in \partial g(\sigma(X)) \cap \mathbb{B}_\delta(\sigma(\overline{W})).$$

Moreover, Proposition 2 also implies that there exists an orthogonal pair $(U, V) \in \mathbb{O}_\mathbb{X}(X) \cap \mathbb{O}_\mathbb{X}(W)$, that is,

$$X = U[\text{Diag}(\sigma(X)) \; 0]V^T \quad \text{and} \quad W = U[\text{Diag}(\sigma(W)) \; 0]V^T.$$

Observe that $(\partial g)^{-1}(\sigma(\overline{W})) = \partial g^*(\sigma(\overline{W}))$ is a nonempty closed convex set [48, Theorem 23.2]. Let $u_X = \Pi_{(\partial g)^{-1}(\sigma(\overline{W}))}(\sigma(X))$, the projection of $\sigma(X)$ onto $(\partial g)^{-1}(\sigma(\overline{W}))$. We deduce from (18) and (20) that

$$\|\sigma(X) - u_X\| = \mathrm{dist}(\sigma(X), (\partial g)^{-1}(\sigma(\overline{W}))) \leq \kappa_0 \|\sigma(\overline{W}) - \sigma(W)\| \leq \kappa_0 \|W - \overline{W}\|.$$

By shrinking $\delta$ if necessary, we know from Lemma 1(i) that there exist $\kappa_1 > 0$ and $(\overline{U}, \overline{V}) \in \mathbb{O}_{\mathbb{X}}(\overline{W})$ such that

$$\|(U, V) - (\overline{U}, \overline{V})\| \leq \kappa_1 \|W - \overline{W}\|. \tag{21}$$

Define $\widehat{X} = \overline{U}[\mathrm{Diag}(u_X)\ 0]\overline{V}^T$ (which is not necessarily a singular value decomposition of $\widehat{X}$). Then $\widehat{X} \in (\partial \theta)^{-1}(\overline{W})$ by Proposition 2. Consequently,

$$\begin{aligned}
\mathrm{dist}(X, (\partial \theta)^{-1}(\overline{W})) &\leq \|X - \widehat{X}\| = \|U[\mathrm{Diag}(\sigma(X))\ 0]V^T - \overline{U}[\mathrm{Diag}(u_X)\ 0]\overline{V}^T\| \\
&\leq \|U[\mathrm{Diag}(\sigma(X))\ 0]V^T - \overline{U}[\mathrm{Diag}(\sigma(X))\ 0]\overline{V}^T\| + \|\sigma(X) - u_X\| \\
&\leq (\|U - \overline{U}\| + \|V - \overline{V}\|)\|\sigma(X)\| + \|\sigma(X) - u_X\| \\
&\leq 2\kappa_1 \|W - \overline{W}\|(\|\sigma(\overline{X})\| + \varepsilon) + \kappa_0 \|W - \overline{W}\| \leq \kappa \|W - \overline{W}\|,
\end{aligned}$$

where $\kappa = 2\kappa_1(\|\sigma(\overline{X})\| + \varepsilon) + \kappa_0$. Since the above inequality is true for any $X \in \mathbb{B}_\varepsilon(\overline{X})$ and any $W \in \partial\theta(X) \cap \mathbb{B}_\delta(\overline{W})$, the inequality (19) then follows. $\square$

By using a similar approach, we can establish an analogous conclusion to Proposition 14 regarding the metric subregularity of the subdifferential of a spectral function associated with a symmetric function in the form of (4).

PROPOSITION 15. *Let $g : \mathbb{R}^n \to (-\infty, +\infty]$ be a proper closed convex and symmetric function. Let $\theta : \mathbb{S}^n \to (-\infty, +\infty]$ be the spectral function associated with $g$ as in (4). Let $(\overline{X}, \overline{W}) \in \mathrm{gph}\,\partial\theta$. If the subdifferential mapping $\partial g$ is metrically subregular at $\lambda(\overline{X})$ for $\lambda(\overline{W})$, then the subdifferential mapping $\partial\theta$ is metrically subregular at $\overline{X}$ for $\overline{W}$.*

*Proof.* One can derive the proof by replacing the singular value decompositions in the proof of Proposition 14 with the eigenvalue decompositions, and then applying Proposition 3 and Lemma 1(ii) accordingly instead of Proposition 2 and Lemma 1(i). For brevity, we omit the details here. $\square$

REMARK 3. *Recall that a set-valued mapping is said to be piecewise polyhedral if its graph is the union of finitely many polyhedral sets. In his Ph.D thesis, Sun [56] shows that a proper closed convex function $p$ is piecewise linear-quadratic if and only if $\partial p$ is piecewise polyhedral, which is also equivalent to $p^*$ being piecewise linear-quadratic (see also [51, Theorem 11.14, Proposition 12.30]). Moreover, a fundamental result of Robinson about the local upper Lipschitz continuous property of the polyhedral mapping [47] implies that for any convex piecewise linear-quadratic function $p$ and any point $(x, w) \in \mathrm{gph}\,\partial p$, $(\partial p)^{-1} = \partial p^*$ is metrically subregular at $x$ for $w$. Thus, the subdifferential of the spectral function $\theta$ is metrically subregular at any $X \in \mathrm{dom}\,\theta$ for any $W \in \partial\theta(X)$ if the underlying function $g$ is a convex piecewise linear-quadratic function. Particular examples of such $g$ include the sum of $k$ largest absolute components of a given vector and the indicator function over a convex polyhedral set. Thus, Proposition 14 and Proposition 15 cover the results in [64, Proposition 11] and [11, Theorem 2.4] for the metric subregularity of the subdifferentials of the nuclear norm function and of the indicator function over the PSD cone, respectively; see also [12, Proposition 3.3] for the metric subregularity of the latter subdifferentials.*

REMARK 4. *For a possibly nonconvex function $p : \mathbb{R}^n \to (-\infty, +\infty]$, one can define its subdifferential mapping $\partial p$ as in [5, (2.227) and (2.229)]. For any $x \in \mathbb{R}^n$, $\partial p(x)$ defined in this way is always a closed and convex set. In fact, all the assertions in Proposition 2 and Proposition 3 are originally established under the more general nonconvex settings in [33, 35]. Hence, one can extend the results in Proposition 14 and Proposition 15 to nonconvex functions $g$ and $\theta$ with no difficulty.*

Now we come back to discuss the second type of sufficient conditions for the quadratic growth conditions of problem (1) and problem (2). It is a well known fact that if the function $h$ is strongly convex on any convex compact subset of $\operatorname{dom} h$, the value $\mathcal{F}\overline{X}$ is invariant over $\overline{X} \in \Omega_P$ (see, e.g., [41]). Hence, the following values are well-defined for any $\overline{X} \in \Omega_P$:

$$\bar{\xi} := \mathcal{F}\overline{X}, \quad \bar{\eta} := \mathcal{F}^* \nabla h(\mathcal{F}\overline{X}) + C. \tag{22}$$

We denote the set $\mathcal{V}_P := \{X \in \mathbb{X} \mid \mathcal{F}X = \bar{\xi}\}$ and two mappings $\mathcal{G}_P^1, \mathcal{G}_P^2 : \mathbb{R}^e \to \mathbb{X}$ as

$$\mathcal{G}_P^1(y) := (\partial \theta)^{-1}(\mathcal{A}^* y - \bar{\eta}), \quad \mathcal{G}_P^2(y) := \{X \in \mathbb{X} \mid 0 \in \mathcal{A}X - b + \mathcal{N}_{\mathcal{Q}^*}(y)\}, \quad y \in \mathbb{R}^e.$$

By the KKT condition (9), it is easy to obtain that

$$\Omega_P = \mathcal{V}_P \cap \mathcal{G}_P^1(\bar{y}) \cap \mathcal{G}_P^2(\bar{y}) \quad \forall \, \bar{y} \in \mathcal{M}_P(\overline{X}).$$

Similarly, if $h^*$ is assumed to be strongly convex in the dual problem (2), the value $\bar{w}$ is invariant over all $\overline{Z} = (\bar{y}, \bar{w}, \overline{S}) \in \Omega_D$. We denote the set $\mathcal{V}_D := \{Z \in \mathbb{Z} \mid \mathcal{A}^* y + \mathcal{F}^* w + S = C\}$ and three mappings $\mathcal{G}_D^1, \mathcal{G}_D^2, \mathcal{G}_D^3 : \mathbb{X} \to \mathbb{Z}$ as

$$\mathcal{G}_D^1(X) := \{Z \in \mathbb{Z} \mid 0 \in \mathcal{A}X - b + \mathcal{N}_{\mathcal{Q}^*}(y)\}, \quad \mathcal{G}_D^2(X) := \{Z \in \mathbb{Z} \mid w = \bar{w}\},$$
$$\mathcal{G}_D^3(X) := \{Z \in \mathbb{Z} \mid 0 \in S + \partial \theta(X)\}, \quad X \in \mathbb{X}.$$

By the KKT condition (8), we have

$$\Omega_D = \mathcal{V}_D \cap \mathcal{G}_D^1(\overline{X}) \cap \mathcal{G}_D^2(\overline{X}) \cap \mathcal{G}_D^3(\overline{X}) \quad \forall \, \overline{X} \in \mathcal{M}_D(\overline{Z}).$$

The following theorem presents conditions for ensuring the quadratic growth conditions for problem (1), under the metric subregularity of the subdifferential of the (absolutely) symmetric function $g$.

THEOREM 16. *Let $\theta$ be a spectral function in the form of (3) (or (4)). Let $\overline{X}$ be an optimal solution for problem (1) with $\mathcal{M}_P(\overline{X}) \neq \emptyset$. Assume that the following three conditions hold:*
*(a) the function $h$ is strongly convex on any compact subset of $\operatorname{dom} h$;*
*(b) for any $\bar{v} \in \partial g(\sigma(\overline{X}))$ (or $\partial g(\lambda(\overline{X}))$), the mapping $\partial g$ is metrically subregular at $\sigma(\overline{X})$ (or $\lambda(\overline{X})$) for $\bar{v}$;*
*(c) the collection of sets $\{\mathcal{V}_P, \mathcal{G}_P^1(\bar{y}), \mathcal{G}_P^2(\bar{y})\}$ is boundedly linearly regular for some $\bar{y} \in \mathcal{M}_P(\overline{X})$.*
*Then the quadratic growth condition (10) holds at $\overline{X}$ for problem (1).*

*Proof.* Applying the results of Proposition 14 to the function $\theta = g \circ \sigma$ or Proposition 15 to the function $\theta = g \circ \lambda$, we see that for any $\overline{S} \in \partial \theta(\overline{X})$, the mapping $\partial \theta$ is metrically subregular at $\overline{X}$ for $\overline{S}$ under condition (b). Then the conclusion follows from [12, Theorem 3.1] directly. □

The remaining issue is how to verify assumption (c) in Theorem 16. By invoking Proposition 8, we arrive at the following results.

PROPOSITION 17. *The collection of sets $\{\mathcal{V}_P, \mathcal{G}_P^1(\bar{y}), \mathcal{G}_P^2(\bar{y})\}$ is boundedly linearly regular for some $\bar{y} \in \mathcal{M}_P(\overline{X})$ under one of the following two conditions:*

*(a) $\mathcal{G}_P^1(\bar{y})$ is a polyhedral set;*

*(b) there exists $\widehat{X} \in \Omega_P$ such that $\widehat{X} \in \text{ri}\left(\mathcal{G}_P^1(\bar{y})\right)$.*

*Proof.* By Proposition 8 and the facts that $\mathcal{V}_P$ and $\mathcal{G}_P^2(\bar{y})$ are polyhedral sets, we directly get the conclusion under assumption (a). If assumption (b) holds, then $\mathcal{V}_P \cap \text{ri}\left(\mathcal{G}_P^1(\bar{y})\right) \cap \mathcal{G}_P^1(\bar{y}) \neq \emptyset$, and the conclusion also follows. □

Below we make several comments regarding the assumptions imposed in Proposition 17. Firstly, let us illustrate condition (a) by several examples.

- Example 1: $\theta(X) = \|X\|_2$ for any $X \in \mathbb{R}^n$ (the vector 2-norm). Then $\partial \theta^{-1}(X)$ is a polyhedral set for any $X \in \mathbb{R}^m$ [58, 64] and the condition (a) automatically holds.
- Example 2: $\theta(X) = \delta_{\mathbb{S}_+^n}(X)$ for any $X \in \mathbb{S}^n$. Then condition (a) is equivalent to $\text{rank}(\mathcal{A}^*\bar{y} - \bar{\eta}) \geq n - 1$ [12, Proposition 3.2].
- Example 3: $\theta(X) = \|X\|_{(k)}$ for any $X \in \mathbb{X}$. Then in order for $\mathcal{G}_P^1(\bar{y})$ to be a polyhedral set, we either have that $\|\mathcal{A}^*\bar{y} - \bar{\eta}\|_* < k$ and $\sigma_2(\mathcal{A}^*\bar{y} - \bar{\eta}) < 1$, or $\|\mathcal{A}^*\bar{y} - \bar{\eta}\|_* = k$, $\sigma_2(\mathcal{A}^*\bar{y} - \bar{\eta}) < 1$ and $\sigma_m(\mathcal{A}^*\bar{y} - \bar{\eta}) > 0$. This result can be obtained via the characterization of $\partial \|\cdot\|_{(k)}$ in [61, 44].

Condition (b) in Proposition 17 can be viewed as the strict complementarity condition for the generalized equation $0 \in -X + \partial\theta(\mathcal{A}^*y - \bar{\eta})$ at $(\overline{X}, \bar{y})$. In particular, one can find characterizations of this property for $\theta(\cdot) = \delta_{\mathbb{S}_+^n}(\cdot)$ in [5, Example 4.79] and for $\theta(\cdot) = \|\cdot\|_{(k)}$ in [14]. It is worth to mention that in order to satisfy assumption (c) in Theorem 16, it suffices for problem (1) to admit a KKT point $(\widehat{X}, \bar{y})$ possessing the partial strict complementarity condition with respect to the spectral function $\theta$. The $\widehat{X}$ here can be different from the reference point $\overline{X}$.

We can also establish sufficient conditions for the quadratic growth condition (11) of the dual problem (2).

THEOREM 18. *Let $\theta$ be a spectral function in the form of (3) (or (4)). Let $\overline{Z}$ be an optimal solution for problem (2) with $\mathcal{M}_D(\overline{Z}) \neq \emptyset$. Assume that the following three conditions hold:*

*(a) the function $h^*$ is continuously differentiable on $\text{dom}\, h^*$, which is assumed to be a nonempty open convex set, and $h^*$ is also strongly convex on any compact subset of $\text{dom}\, h^*$;*

*(b) for any $\bar{w} \in \partial g^*(\sigma(\overline{S}))$ (or $\partial g^*(\lambda(\overline{S}))$), the mapping $\partial g^*$ is metrically subregular at $\sigma(\overline{S})$ (or $\lambda(\overline{S})$) for $\bar{w}$;*

*(c) the collection of sets $\{\mathcal{V}_D, \mathcal{G}_D^1(\overline{X}), \mathcal{G}_D^2(\overline{X}), \mathcal{G}_D^3(\overline{X})\}$ is boundedly linearly regular for some $\overline{X} \in \mathcal{M}_D(\overline{Z})$.*

*Then the quadratic growth condition (11) holds at $\overline{Z}$ for problem (2).*

Before proceeding to the next section, we mention that when $\mathcal{A}$, $b$ and $\mathcal{Q}$ are vacant in problem (1) and $h$ is taken to be a least square function $\frac{1}{2}\|\cdot - d\|^2$ for some given $d \in \mathbb{R}^d$, the problem (1) reduces to

$$\text{(23)} \qquad \min \frac{1}{2}\|\mathcal{F}X - d\|^2 + \langle C, X \rangle + \theta(X),$$

and the corresponding dual is given by

$$\text{(24)} \qquad \min \frac{1}{2}\|w - d\|^2 + \theta^*(-S), \quad \text{s.t.} \quad \mathcal{F}^*w + S = C.$$

Obviously problem (24) has a unique optimal solution. Denote it as $(\bar{w}, \overline{S})$. Suppose that $\theta$ is given by (3) and $g^*$ is $\mathcal{C}^2$ cone reducible at $-\sigma(\overline{S})$, then the quadratic growth condition for the dual problem (24) always holds at $(\bar{w}, \overline{S})$ by Theorem 13. However, certain conditions need to be imposed to the primal form (23) for deriving its quadratic growth condition, like the one given in Theorem 12 or Theorem 16. Examples for the failure of the quadratic growth condition of problem (23) without any additional assumption can be found in [64], where $\theta$ is taken to be the nuclear norm function.

**4. Fast convergence rates of the ALM.** In recent years, lots of progress has been achieved in solving large scale semidefinite programming (SDP) and non-symmetric matrix optimization problems regularized by the nuclear norm and the spectral norm [63, 32, 10, 62, 38]. The central idea of these works is to combine the ALM and the semismooth matrix-valued function theory. In this section, we demonstrate that under the derived quadratic growth conditions in Theorem 12 and Theorem 16, the iteration sequence generated by the ALM enjoys fast convergence rates. This part of work is motivated by a recent paper [12] on discussing the asymptotic (super)linear convergence rates of the ALM for solving the composite SDP.

Recall that we denote the variable $Z = (y, w, S)$ for any $y \in \mathbb{R}^e$, $w \in \mathbb{R}^d$ and $S \in \mathbb{X}$, and the space $\mathbb{Z} = \mathbb{R}^e \times \mathbb{R}^d \times \mathbb{X}$. Let $c > 0$ be a positive parameter. For any $Z \in \mathbb{Z}$ and $X \in \mathbb{X}$, the augmented Lagrangian function associated with problem (2) is given by

$$(25) \qquad L_c(Z, X) := l(Z, X) + (c/2)\|\mathcal{A}^*y + \mathcal{F}^*w + S - C\|^2.$$

Let $k \geqslant 0$. Given a sequence of scalars $c_k \uparrow c_\infty \leqslant \infty$ and a starting point $X^0 \in \mathbb{X}$, the $(k+1)$-th iteration of the ALM is given by

$$(26) \qquad \begin{cases} Z^{k+1} \approx \arg\min_{Z \in \mathbb{Z}} \left\{ \zeta_k(Z) := L_{c_k}(Z, X^k) \right\}, \\ X^{k+1} = X^k + c_k(\mathcal{A}^*y^{k+1} + \mathcal{F}^*w^{k+1} + S^{k+1} - C), \end{cases} \quad k \geqslant 0.$$

Our approach for deriving the convergence rates of the ALM is based on the observation by Rockafellar that the ALM is a special implementation of the inexact dual proximal point algorithm (PPA) for solving convex optimization problems [49]. Specifically, let $P_k := (\mathcal{I} + c_k \mathcal{T}_\phi)^{-1}$, where $\mathcal{I}$ is the identity mapping in $\mathbb{X}$ and $\mathcal{T}_\phi$ is the set-valued mapping given in (7). Consider the PPA

$$(27) \qquad X^{k+1} \approx P_k(X^k)$$

to be executed with the following stopping criteria:

(A) $\quad \|X^{k+1} - P_k(X^k)\| \leqslant \varepsilon_k, \quad \varepsilon_k \geqslant 0, \quad \sum_{k=0}^{\infty} \varepsilon_k < \infty,$

(B) $\quad \|X^{k+1} - P_k(X^k)\| \leqslant \eta_k \|X^{k+1} - X^k\|, \quad \eta_k \geqslant 0, \quad \sum_{k=0}^{\infty} \eta_k < \infty.$

Then the relationship between the iteration sequences generated by the ALM in (26) for solving problem (2) and by the PPA in (27) for solving problem (1) is demonstrated in the following proposition. We prove this proposition by adapting the proofs in [49, Proposition 6] to our settings.

PROPOSITION 19. *Given $X^k \in \mathbb{X}$, $Z^{k+1} = (y^{k+1}, w^{k+1}, S^{k+1}) \in \mathbb{Z}$ and a positive parameter $c_k$ for some $k \geqslant 0$. Let $\zeta_k$ and $P_k$ be given by (26) and (27). Denote $X^{k+1} = X^k + c_k(\mathcal{A}^*y^{k+1} + \mathcal{F}^*w^{k+1} + S^{k+1} - C)$. Then*

$$(28) \qquad \|X^{k+1} - P_k(X^k)\|^2/(2c_k) \leqslant \zeta_k(Z^{k+1}) - \inf \zeta_k.$$

*Proof.* Firstly, it is known from [48, Theorem 37.3] that for any $c > 0$ and $X \in \mathbb{X}$,

$$
\begin{aligned}
\inf_Z L_c(Z, X) &= \inf_Z \max_{\widehat{X} \in \mathbb{X}} \{l(Z, \widehat{X}) - 1/(2c)\|\widehat{X} - X\|^2\} \\
&= \max_{\widehat{X} \in \mathbb{X}} \inf_Z \{l(Z, \widehat{X}) - 1/(2c)\|\widehat{X} - X\|^2\} \\
&= \max_{\widehat{X} \in \mathbb{X}} \{\phi(\widehat{X}) - 1/(2c)\|\widehat{X} - X\|^2\},
\end{aligned}
\tag{29}
$$

where $\phi$ is the essential objective function of problem (1) defined by (6). By taking into account the definition of $P_k$ in (27), we see that

$$L_{c_k}(Z^{k+1}, X) \geqslant \inf_Z L_{c_k}(Z, X) \geqslant \phi(P_k(X^k)) - 1/(2c_k)\|P_k(X^k) - X\|^2 \quad \forall\, X \in \mathbb{X}$$

and

$$\inf_Z L_{c_k}(Z, X^k) = \phi(P_k(X^k)) - 1/(2c_k)\|P_k(X^k) - X^k\|^2.$$

In view of the above two inequalities, we know that for any $X \in \mathbb{X}$,

$$
\begin{aligned}
&\zeta_k(Z^{k+1}) - \inf \zeta_k \\
={}& L_{c_k}(Z^{k+1}, X^k) - \inf_Z L_{c_k}(Z, X^k) \\
={}& L_{c_k}(Z^{k+1}, X) - \langle X^{k+1} - X^k, X - X^k \rangle/c_k - \inf_Z L_{c_k}(Z, X^k) \\
\geqslant{}& \phi(P_k(X^k)) - 1/(2c_k)\|P_k(X^k) - X\|^2 - \langle X^{k+1} - X^k, X - X^k \rangle/c_k \\
&- (\phi(P_k(X^k)) - 1/(2c_k)\|P_k(X^k) - X^k\|^2) \\
={}& (2\langle P_k(X^k) - X^{k+1}, X - X^k \rangle - \|X - X^k\|^2)/(2c_k).
\end{aligned}
$$

Then (28) can be established by taking $X = P_k(X^k) - X^{k+1} + X^k$ in the above inequality. □

It follows from the above proposition that the sequence of multiplies $\{X^k\}$ generated by the ALM for solving problem (2), if adopted the following two stopping criteria,

$(A')$ $\zeta_k(Z^{k+1}) - \inf \zeta_k \leqslant \varepsilon_k^2/2c_k$, $\quad \varepsilon_k \geqslant 0$, $\quad \sum_{k=0}^{\infty} \varepsilon_k < \infty$,

$(B')$ $\zeta_k(Z^{k+1}) - \inf \zeta_k \leqslant (\eta_k^2/2c_k)\|X^{k+1} - X^k\|^2$, $\quad \eta_k \geqslant 0$, $\quad \sum_{k=0}^{\infty} \eta_k < \infty$,

can be taken as the iterates generated by the inexact PPA for solving problem (1) with the stopping criteria $(A)$ and $(B)$, respectively. This fact allows us to study the convergence rates of the ALM via the rates of the inexact PPA.

In [50], Rockafellar established the convergence rates of the inexact PPA under the Lipschitz continuity of $\mathcal{T}_\phi^{-1}$ at the origin. This assumption, which automatically requires the uniqueness of the optimal solution of problem (1), has been relaxed by Luque with an error bound type condition [40, (2.1)]. Luque's condition is known to be satisfied if $\mathcal{T}_\phi$ is a polyhedral mapping [47]. When $\mathcal{T}_\phi$ is non-polyhedral, to check Luque's condition may be difficult, especially for the case that $\mathcal{T}_\phi^{-1}(0)$ is unbounded. In [12], this condition is further relaxed by the metric subregularity of the operator $\mathcal{T}_\phi$ at some optimal point for the origin. Based on Proposition 19 and [12, Theorem 4.1], we can prove the global and local (super)linear convergence of the ALM for solving matrix optimization problem (1).

THEOREM 20. *Assume that the optimal solution set $\Omega_P$ to problem (1) is non-empty. Denote $\Psi^*$ be the optimal value of $\Psi$ for problem (2). Let $\{(Z^k, X^k)\}$ be an* *infinite sequence generated by the ALM in (26) with stopping criterion $(A')$, where $Z^k = (y^k, w^k, S^k)$. Then, the whole sequence $\{X^k\}$ is bounded and converges to some $X^\infty \in \Omega_P$, and the sequence $\{Z^k\}$ satisfies for all $k \geqslant 0$,*

$$\text{(30)} \qquad \|\mathcal{A}^* y^{k+1} + \mathcal{F}^* w^{k+1} + S^{k+1} - C\| = c_k^{-1}\|X^{k+1} - X^k\| \to 0,$$

$$\text{(31)} \qquad \Psi(Z^{k+1}) - \Psi^* \leqslant \zeta_k(Z^{k+1}) - \inf \zeta_k + (1/2c_k)(\|X^k\|^2 - \|X^{k+1}\|^2).$$

*Moreover, if problem (2) has a non-empty bounded solution set, then the sequence $\{Z^k\}$ is also bounded, and all of its accumulation points are optimal solutions to problem (2).*

*If the quadratic growth condition (10) holds at $X^\infty$ for the origin with modulus $\kappa_p > 0$, then under criterion $(B')$: there exists $\bar{k} \geqslant 0$ such that for all $k \geqslant \bar{k}$,*

$$\text{(32)} \qquad \operatorname{dist}(X^{k+1}, \Omega_P) \leqslant \theta_k \operatorname{dist}(X^k, \Omega_P),$$

*where*

$$\theta_k = (\mu_k + 2\eta_k)(1 - \eta_k)^{-1} \quad \text{with} \quad \mu_k = 1/\sqrt{1 + c_k^2 \kappa_p^2},$$
$$\theta_k \to \theta_\infty = 1/\sqrt{1 + c_\infty^2 \kappa_p^2} \quad (\theta_\infty = 0 \text{ if } c_\infty = \infty).$$

*In addition, the following inequalities hold regarding the R-(super)linear convergence rate of dual feasibility and dual objective value:*

$$\text{(33a)} \qquad \|\mathcal{A}^* y^{k+1} + \mathcal{F}^* w^{k+1} + S^{k+1} - C\| \leqslant \tau_k^1 \operatorname{dist}(X^k, \Omega_P),$$
$$\text{(33b)} \qquad \Psi(Z^{k+1}) - \Psi^* \leqslant \tau_k^2 \operatorname{dist}(X^k, \Omega_P),$$

*where*

$$\tau_k^1 := c_k^{-1}(1 - \eta_k)^{-1} \to \tau_\infty^1 = 1/c_\infty,$$
$$\tau_k^2 := \tau_k^1(\eta_k^2 \|X^{k+1} - X^k\| + \|X^{k+1}\| + \|X^k\|)/2 \to \tau_\infty^2 = \|X^\infty\|/c_\infty,$$
$$(\tau_\infty^1 = \tau_\infty^2 = 0 \text{ if } c_\infty = \infty).$$

*Proof.* The assertions in the first paragraph of Theorem 20 can be obtained from [49, Theorem 4] and [12, Theorem 4.2]. From Lemma 9 we know that if the quadratic growth condition (10) holds at $X^\infty$ for the origin with modulus $\kappa_p > 0$, then $\mathcal{T}_\phi$ is metrically subregular at $X^\infty$ for the origin with modulus $1/\kappa_p$. Thus, the inequality (32) follows from Proposition 19 and [12, Theorem 4.1]. Finally, one can establish the inequalities (33a) and (33b) by adapting the proofs in [12, Theorem 4.2] to our settings with no difficulty. □

Theorem 20 shows that under the quadratic growth condition of problem (1), the primal iteration sequence $\{X^k\}$ generated by the ALM converges Q-(super)linearly, while the dual feasibility and objective value converge at least R-(super)linearly. Perhaps more importantly, the Q-(super)linear rate constant $\theta_k$ of the sequence $\{X^k\}$ could be much smaller than 1. For example, the rate $\theta_k$ shall be around $\sqrt{2}/2$ when the penalty parameter $c_k$ approaches $1/\kappa_p$. This distinctive feature makes the ALM very attractive. Numerical experiments conducted in the next section will demonstrate this point.

**5. Numerical experiments.** In this section, we conduct numerical experiments for the ALM on solving the fastest mixing Markov chain (FMMC) problem [6, 7] in order to illustrate its fast convergence rates. Let $\mathcal{G} = (\mathcal{V}, \mathcal{E})$ be a connected graph with vertex set $\mathcal{V} = \{1, \ldots, n\}$ and edge set $\mathcal{E} \subseteq \mathcal{V} \times \mathcal{V}$. Label the edges by $l = 1, \ldots, d$. For any $y \in \Re^d$, denote $\mathcal{A}^* y = \sum_{l=1}^{d} y_l E^{(l)}$, where if the edge $l$ connects two vertices $i$ and $j$ ($i \neq j$), then $E^{(l)}_{ij} = E^{(l)}_{ji} = 1$, $E^{(l)}_{ii} = E^{(l)}_{jj} = -1$ and all other entries of $E^{(l)}$ are zero. Let $B \in \mathbb{R}^{n \times d}$ be the vertex-edge incidence matrix defined by

$$B_{il} = \begin{cases} 1 & \text{if edge } l \text{ incident to vertx } i, \\ 0 & \text{otherwise.} \end{cases} \quad i \in \{1, \ldots, n\} \text{ and } j \in \{1, \ldots, d\}.$$

The corresponding FMMC problem can be equivalent written in terms of the optimization variables $y \in \mathbb{R}^d$, $z \in \mathbb{R}^{n+d}$ and $P \in \mathbb{S}^n$ as

$$(34) \qquad \begin{aligned} \min \quad & \|P\|_{(2)} + \delta_{\mathbb{R}_+^{n+d}}(z) \\ \text{s.t.} \quad & -\mathcal{A}^* y + P = I_n, \quad \overline{B} y + z = \bar{e}, \end{aligned}$$

where $\|\cdot\|_{(2)}$ is the Ky Fan 2-norm in $\mathbb{S}^n$, $\overline{B} = \begin{bmatrix} I_d \\ -B \end{bmatrix} \in \mathbb{R}^{(n+d) \times d}$ and $\bar{e} = (0, e_n)^T \in \mathbb{R}^{n+d}$. One can easily see that problem (34) is in the form of (2) by taking the variable $S$ in (2) as $(P, z)$, the set $\mathcal{Q}^* = \{0\} \subseteq \mathbb{R}^d$, the function $\theta^* = \|\cdot\|_{(2)} \times \delta_{\mathbb{R}_+^{n+d}}$, the constant $b = 0$ and the linear operator $\mathcal{F}^* = 0$. Therefore, we can adopt the ALM discussed in the last section for solving problem (34).

**5.1. Solving the subproblems of the ALM.** Given the penalty parameter $c > 0$, the augmented Lagrangian function of problem (34) takes the form of

$$\begin{aligned} L_c(P, z, y, X_1, X_2) = & \|P\|_{(2)} + \delta_{\mathbb{R}_+^{n+d}}(z) + \langle X_1, -\mathcal{A}^* y + P - I_n \rangle + \langle X_2, \overline{B} y + z - \bar{e} \rangle \\ & + \tfrac{c}{2} \| -\mathcal{A}^* y + P - I_n \|^2 + \tfrac{c}{2} \| \overline{B} y + z - \bar{e} \|^2, \end{aligned}$$

where $(P, z, y, X_1, X_2) \in \mathbb{S}^n \times \mathbb{R}^{n+d} \times \mathbb{R}^d \times \mathbb{S}^n \times \mathbb{R}^{n+d}$. As discussed in Section 4, the ALM for solving problem (34) takes the following iteration:

$$\begin{cases} (P^{k+1}, z^{k+1}, y^{k+1}) \approx \arg\min\{\zeta_k(P, z, y) := L_{c_k}(P, z, y, X_1^k, X_2^k)\}, \\ (X_1^{k+1}, X_2^{k+1}) = (X_1^k, X_2^k) + c_k(-\mathcal{A}^* y^{k+1} + P^{k+1} - I_n, \overline{B} y^{k+1} + z^{k+1} - \bar{e}) \end{cases}$$

for a sequence of scalars $c_k \uparrow c_\infty \leq \infty$ with $k \geq 0$. Obviously, the major computational cost of the above framework comes from obtaining approximate solutions of the subproblems.

We shall adopt the semismooth Newton-CG method to solve the subproblems as in [63, 10, 62, 38]. Denote, for any $k \geq 0$ and $y \in \mathbb{R}^d$,

$$P_k(y) = \text{Prox}_{\|\cdot\|_{(2)}/c_k}(\mathcal{A}^* y + I_n - X_1^k/c_k), \quad z_k(y) = \Pi_{\mathbb{R}_+^{n+d}}(-\overline{B} y + \bar{e} - X_2^k/c_k).$$

Then $(\widetilde{P}, \tilde{z}, \tilde{y}) \in \arg\min \zeta_k(P, z, y)$ if and only if

$$\tilde{y} \in \arg\min\{\xi_k(y) := \zeta_k(P_k(y), z_k(y), y)\}, \quad \widetilde{P} = P_k(\tilde{y}), \quad \tilde{z} = z_k(\tilde{y}).$$

The function $\xi_k$ is continuously differentiable with the gradient given by (see, e.g., [30]):

$$\nabla \xi_k(y) = \mathcal{A} \Pi_{\mathbb{B}^*}(c_k(\mathcal{A}^* y + I_n) - X_1^k) + \overline{B}^T \Pi_{\mathbb{R}_+^{n+d}}(c_k(\overline{B}\bar{y} - \bar{e}) + X_2^k),$$

where $\mathbb{B}^* \subseteq \mathbb{S}^n$ is the unit ball associated with the dual norm $\|\cdot\|_{(2)}^*$, i.e.,

$$\mathbb{B}^* := \{X \in \mathbb{S}^n \mid \|X\|_{(2)}^* \leq 1\} = \{X \in \mathbb{S}^n \mid \|X\|_2 \leq 1, \ \|X\|_* \leq 2\}.$$

Since $\Pi_{\mathbb{B}^*}(\cdot)$ is a matrix spectral operator defined in [15], its strong semismoothness then follows from [20, Proposition 7.4.7] and [13, Theorem 3.12]. Thus, $\nabla \xi_k$ is strongly semismooth as it is the composition of strongly semismooth functions [21, 54]. This allows us to develop a semismooth Newton-CG method for solving the nonlinear equation $\nabla \xi_k(\cdot) = 0$. For the implementation and the convergence analysis of the semismooth Newton-CG method, the interested reader may refer to [63, 10] for details.

**5.2. Easy-to-verify stopping criteria.** The stopping criteria $(A')$ and $(B')$ can be difficult to execute due to the involvement of the unknown values $\inf \zeta_k$. For the practical implementation and numerical considerations, here we provide easy ways to verify these two criteria in this subsection. This part is inspired by the work done in [12] on the implementation of the ALM for solving the convex SDP problems.

Note that the dual of problem (34) is given by

(35)
$$\begin{aligned}\max \quad & \langle I, X_1 \rangle - \langle \bar{e}, X_2 \rangle \\ \text{s.t.} \quad & \mathcal{A}X_1 - \overline{B}^T X_2 = 0, \ X_1 \in \mathbb{B}^*, \ X_2 \in \mathbb{R}_+^{n+d}.\end{aligned}$$

For any $k \geq 0$, denote the concave function $f_k : \mathbb{S}^n \times \mathbb{R}^{n+d} \to (-\infty, +\infty)$ as

$$f_k(X) := \langle I, X_1 \rangle - \langle \bar{e}, X_2 \rangle - \frac{1}{2c_k}(\|X_1 - X_1^k\|^2 + \|X_2 - X_2^k\|^2),$$

where $X := (X_1, X_2) \in \mathbb{S}^n \times \mathbb{R}^{n+d}$.

ASSUMPTION 1. *Assume that there exists $\widehat{X} = (\widehat{X}_1, \widehat{X}_2)$ such that*

$$\mathcal{A}\widehat{X}_1 - \overline{B}^T \widehat{X}_2 = 0, \ \widehat{X}_1 \in \text{int}(\mathbb{B}^*), \ \widehat{X}_2 \in \mathbb{R}_+^{n+d}.$$

Let $F$ be the feasible set of problem (35), i.e.,

$$F := \{(X_1, X_2) \in \mathbb{S}^n \times \mathbb{R}^{n+d} \mid \mathcal{A}X_1 - \overline{B}^T X_2 = 0, \ X_1 \in \mathbb{B}^*, \ X_2 \in \mathbb{R}_+^{n+d}\}.$$

It is known from [4, Theorem 7] that under Assumption 1, there exists a positive constant $\bar{\mu}$ such that

(36) $$\|X - \Pi_F(X)\| \leq \bar{\mu}(1 + \|X\|) r(X) \quad \forall X \in \mathbb{S}^n \times \mathbb{R}^{n+d},$$

where

$$r(X) := \|\mathcal{A}X_1 - \overline{B}^T X_2\| + \|X_1 - \Pi_{\mathbb{B}^*}(X_1)\| + \|X_2 - \Pi_{\mathbb{R}_+^{n+d}}(X_2)\|$$

for any $X = (X_1, X_2) \in \mathbb{S}^n \times \mathbb{R}^{n+d}$. Let $\{\varepsilon_k\}$ and $\{\eta_k\}$ be two positive summable sequences for $k \geq 0$. Consider the following two stopping criteria:

$(A'')$ $\begin{cases} \zeta_k(\xi^{k+1}) - f_k(X^{k+1}) \leq \varepsilon_k^2/(2c_k), \\ (1 + \|X^{k+1}\|) r(X^{k+1}) \leq \min\left\{1, \dfrac{\varepsilon_k}{2c_k \|\nabla f_k(X^{k+1})\|}\right\} \varepsilon_k, \end{cases}$

$(B'')$ $\begin{cases} \zeta_k(\xi^{k+1}) - f_k(X^{k+1}) \leq \eta_k^2 \|X^{k+1} - X^k\|^2/(2c_k), \\ (1 + \|X^{k+1}\|) r(X^{k+1}) \leq \min\left\{1, \dfrac{\eta_k}{2c_k \|\nabla f_k(X^{k+1})\|}\right\} \|X^{k+1} - X^k\| \eta_k. \end{cases}$

The following proposition can be derived in the same fashion as in [12, Proposition 4.3].

PROPOSITION 21. *Suppose that Assumption 1 holds. Let $\bar{\mu}$ be given by (36). Suppose that for some $k \geq 0$, $(X_1^k, X_2^k) \in \mathbb{S}^n \times \mathbb{R}^{n+d}$ is not an optimal solution to problem (35). Let $\{\xi^{k,j}\}_{j \geq 0} := \{(P^{k,j}, z^{k,j}, y^{k,j})\}_{j \geq 0}$ be any sequence such that $\zeta_k(\xi^{k,j}) \to \inf \zeta_k$. For any $j \geq 0$, let*

$$X^{k,j} := (X_1^{k,j}, X_2^{k,j}) = (X_1^k, X_2^k) + c_k(-\mathcal{A}^* y^{k,j} + P^{k,j} - I_n, \overline{B} y^{k,j} + z^{k,j} - \bar{e}).$$

*Then there exist $j_A \geq 0$ and $j_B \geq 0$ such that $(A'')$ and $(B'')$ are satisfied by some $(\xi^{k,j_A}, X^{k,j_A})$ and $(\xi^{k,j_B}, X^{k,j_B})$, respectively. Moreover,*

$$(37) \qquad \zeta_k(\xi^{k,j_A}) - \inf \zeta_k \leq \left(1 + \bar{\mu} + \frac{1}{2}\bar{\mu}^2\right) \varepsilon_k^2/(2c_k)$$

*and*

$$(38) \qquad \zeta_k(\xi^{k,j_B}) - \inf \zeta_k \leq \left(1 + \bar{\mu} + \frac{1}{2}\bar{\mu}^2\right) \eta_k^2 \|X^{k,j_B} - X^k\|^2/(2c_k).$$

Proposition 21 indicates that stopping criteria $(A')$ and $(B')$ can be executed via the easy-to-verify criteria $(A'')$ and $(B'')$, respectively.

**5.3. Numerical examples.** To illustrate that the fast convergence rates of the ALM are fundamentally different from the linear convergence rates of first order methods, we also solve the problems by using an alternating direction method of multipliers (ADMM) of Glowinski and Marroco [27] and Gabay and Mercier [23]. Let $(y^0, X_1^0, X_2^0) \in \mathbb{R}^d \times \mathbb{B}^* \times \mathbb{R}_+^{n+d}$ be an initial point. The $(k+1)$-th iteration $(k \geq 0)$ of the ADMM is given by

$$\begin{cases} (P^{k+1}, z^{k+1}) = \arg\min \{L_c(P, z, y^k, X_1^k, X_2^k)\}, \\ y^{k+1} = \arg\min \{L_c(P^{k+1}, z^{k+1}, y, X_1^k, X_2^k)\}, \\ (X_1^{k+1}, X_2^{k+1}) = (X_1^k, X_2^k) + \tau c(-\mathcal{A}^* y^{k+1} + P^{k+1} - I_n, \overline{B} y^{k+1} + z^{k+1} - \bar{e}), \end{cases}$$

where $\tau \in (0, \frac{1+\sqrt{5}}{2})$ is the step-length.

In the numerical experiments, the accuracy of an approximate KKT solution $(P, z, y, X_1, X_2)$ is measured by the relative residual of the primal infeasibility, the dual infeasibility and the duality gap:

$$\eta = \max\{\eta_p, \eta_d, \eta_{gap}\},$$

where

$$\begin{cases} \eta_p = \max\left\{\frac{\|-\mathcal{A}^* y + P - I_n\| + \|\overline{B} y + z - \bar{e}\|}{1 + \|I_n\| + \|\bar{e}\|}, \frac{\|\max\{-z, 0\}\|}{1 + \|z\|}\right\}, \\ \eta_d = \max\left\{\frac{\|\mathcal{A} X_1 - \overline{B}^T X_2\|}{1 + \|\mathcal{A}\| + \|\overline{B}\|}, \frac{\|X_1 - \Pi_{\mathbb{B}^*}(X_1)\|}{1 + \|X_1\|}, \frac{\|\max\{-X_2, 0\}\|}{1 + \|X_2\|}\right\}, \\ \eta_{gap} = \frac{|\|P\|_{(2)} - \langle I, X_1 \rangle + \langle \bar{e}, X_2 \rangle|}{1 + \|P\|_{(2)} + |\langle I, X_1 \rangle - \langle \bar{e}, X_2 \rangle|}. \end{cases}$$

We shall terminate both the ADMM and the ALM when $\eta < 10^{-6}$ with the maximum number of iterations set at 25,000 for the ADMM and 100 for the ALM, respectively. All our numerical results are obtained from a workstation running on a 64-bit Windows Operating System having 12 cores with Intel Xeon E5-2680 processors at 2.50GHz and 128 GB memory. We have implemented the algorithms in Matlab version 9.0. The tested graphs are taken from the SuiteSparse Matrix Collection http://www.cise.ufl.edu/research/sparse/matrices/.

In our implementation of the ADMM, the step-length $\tau$ is set to be 1.618. For the implementation of the ALM, we first run the ADMM to generate a starting point: if the point generated by the ADMM satisfies $\eta < 10^{-4}$ or if the ADMM reaches the number of 200 iterations, the algorithm is switched to the ALM iterations. The stopping criterion $(B'')$ discussed in Section 5.2 is adopted in order to achieve the fast convergence rates.

Table 1: The performance of the ALM and the ADMM for solving the FMMC problems. The computational time is in the format of "hours:minutes:seconds".

|  |  | iteration | $\eta$ | time |
|---|---|---|---|---|
| problem | $d; n$ | ALM \| ADMM | ALM \| ADMM | ALM \| ADMM |
| cage | 2562 ; 366 | 6;6;200 \| 1925 | 0.0-7 \| 9.1-7 | 05 \| 37 |
| G3 | 19176 ; 800 | 32;57;88 \| 599 | 3.0-7 \| 8.3-7 | 1:37 \| 1:21 |
| G6 | 9665 ; 800 | 30;44;145 \| 989 | 8.5-7 \| 9.6-7 | 1:09 \| 1:52 |
| G15 | 4661 ; 800 | 31;51;200 \| 6122 | 3.4-8 \| 7.7-7 | 1:05 \| 11:27 |
| G46 | 9990 ; 1000 | 30;44;134 \| 1619 | 5.6-7 \| 9.9-7 | 1:34 \| 5:33 |
| G54 | 5916 ; 1000 | 22;62;200 \| 8928 | 7.7-7 \| 9.9-7 | 2:23 \| 27:16 |
| G43 | 9990 ; 1000 | 24;96;90 \| 2073 | 2.9-7 \| 9.3-7 | 2:37 \| 6:03 |
| delaunayn10 | 3056 ; 1024 | 61;359;200 \| 25000 | 6.8-9 \| 7.2-5 | 10:25 \| 1:04:29 |
| G22 | 19990 ; 2000 | 31;46;56 \| 2918 | 2.3-8 \| 9.9-7 | 5:57 \| 41:31 |
| G24 | 19990 ; 2000 | 41;296;200 \| 6808 | 2.7-7 \| 9.9-7 | 53:57 \| 1:38:37 |
| G26 | 19990 ; 2000 | 29;87;200 \| 2954 | 1.4-7 \| 9.9-7 | 16:15 \| 42:18 |
| minnesota | 3303 ; 2642 | 25;24;123 \| 258 | 0.0-10 \| 9.1-7 | 6:27 \| 6:29 |
| G48 | 6000 ; 3000 | 40;79;200 \| 9470 | 9.2-7 \| 9.2-7 | 19:39 \| 4:40:08 |
| G49 | 6000 ; 3000 | 25;38;200 \| 7488 | 6.6-7 \| 8.5-7 | 10:58 \| 3:36:58 |
| G50 | 6000 ; 3000 | 26;42;74 \| 6370 | 4.8-8 \| 7.8-7 | 9:36 \| 2:59:34 |
| USpowerGrid | 6594 ; 4941 | 27;120;200 \| 25000 | 1.4-7 \| 1.0-5 | 3:11:22 \| 56:47:18 |

In Table 1, we report the numerical results obtained by the ALM and the ADMM in solving various instances of the FMMC problem (34). The three numbers in the iteration column of the ALM represent the number of outer ALM iterations, the number of inner semismooth Newton-CG iterations and the number of the iterations for generating a starting point, respectively. The fast convergence rates of the ALM can be observed clearly from the very small number of outer iterations. We can also observe from Table 1 that for very easy problems such as *cage*, *G3* and *G6*, the ADMM can also be efficient. This is because the ADMM may also converge linearly under a KKT-type error bound condition [29]. However, this error bound condition is more restrictive than either the primal or the dual second order growth condition for matrix optimization problems. For such examples, see, e.g., [12, Example 1]. Moreover, even if the ADMM possesses the linear convergence, the linear rate constant can be very close to 1 for difficult instances, regardless of the value of the penalty parameter $c$ [29]. These two inherent drawbacks of the ADMM may explain why it does not perform as good as the ALM for the difficult problems listed in Table 1. In particular, for the instance *USpowerGrid*, the ALM solves the problem in about 3 hours, while the ADMM is not able to achieve the required accuracy within 25,000 iterations after more than 56 hours.

**6. Conclusion.** The quadratic growth conditions play a central role in the study of optimization problems, both for perturbation theory and for convergence analysis of optimization algorithms. In this paper, we establish two types of sufficient conditions for ensuring the quadratic growth conditions for a wide class of convex matrix optimization problems associated with spectral functions. One type is based on the "no-gap" second order sufficient conditions of matrix optimization problems, under the $\mathcal{C}^2$-cone reducibility of spectral functions. The other type is through the bounded linear regularity of a collection of sets, under the metric subregularity of subdifferentials of spectral functions. Moreover, we show that these two variational properties of spectral functions, namely, the $\mathcal{C}^2$-cone reducibility and the metric subregularity of the subdifferentials, can be verified via the corresponding properties of underlying (absolutely) symmetric functions. Finally, the quadratic growth conditions are applied to conduct convergence rates analysis of the ALM for solving convex matrix optimization problems. Nevertheless, there remain many issues that require further investigation. These include the sufficient conditions for ensuring the calmness of the KKT solution mappings of matrix optimization problems and the efficient algorithms for solving the subproblems in the ALM when the problems are degenerate.

**Acknowledgements** We would like to thank the two anonymous referees for their helpful comments on improving the quality of this paper. Thanks also go to Professor Defeng Sun, Professor Kim-Chuan Toh and Dr. Xudong Li for many helpful discussions on this work.